\documentclass[12pt,a4paper]{article}
\usepackage[tbtags]{amsmath}    
\usepackage{amsfonts,amssymb,amsthm,euscript,makeidx,pifont,boxedminipage,multicol,fullpage}
\usepackage{color}

\def\be{\begin{eqnarray}}
\def\ee{\end{eqnarray}}

\def\b*{\begin{eqnarray*}}
\def\e*{\end{eqnarray*}}

\newtheorem{Theorem}{Theorem}[section]
\newtheorem{Lemma}[Theorem]{Lemma}
\newtheorem{Proposition}[Theorem]{Proposition}

\newtheorem{Remark}[Theorem]{Remark}
\newtheorem{Example}[Theorem]{Example}

\newtheorem{Assumption}[Theorem]{Assumption}

\newcommand{\rmi}{{\rm (i)$\>\>$}}
\newcommand{\rmii}{{\rm (ii)$\>\>$}}
\newcommand{\rmiii}{{\rm (iii)$\>\>$}}
\newcommand{\rmiv}{{\rm (iv)$\>\>$}}
\newcommand{\rmv}{{\rm (v)$\>\>$}}
\newcommand{\rmvi}{{\rm (vi)$\>\>$}}

\def \E{\mathbb{E}}
\def \F{\mathbb{F}}
\def \G{\mathbb{G}}
\def \L{\mathbb{L}}
\def \P{\mathbb{P}}
\def \Q{\mathbb{Q}}
\def \R{\mathbb{R}}

\def \M{\mathbf{M}}
\def \N{\mathbb{N}}

\def\Ac{{\cal A}}
\def\Bc{{\cal B}}
\def\Cc{{\cal C}}
\def\Dc{{\cal D}}

\def\Fc{{\cal F}}
\def\Gc{{\cal G}}
\def\Hc{{\cal H}}

\def\Mc{{\cal M}}
\def\Nc{{\cal N}}
\def\Oc{{\cal O}}
\def\Pc{{\cal P}}

\def\Sc{{\cal S}}
\def\Tc{{\cal T}}

\def\Wc{{\cal W}}

\def\Pt{{\tilde P}}

\def \Om{\Omega}
\def \om{\omega}
\def \Omb{\overline{\Om}}
\def \omb{\bar{\om}}

\def \eps{\varepsilon}

\def \0{\mathbf{0}}
\def \1{\mathbf{1}}
\def \x{\times}

\def\Xcr{\langle X \rangle}

\def \Fcb{\overline{{\cal F}}}
\def \Fbb{\overline{\F}}
\def \Fct{\tilde{{\cal F}}}
\def \Ft{\tilde{\F}}
\def \Pcb{\overline{\Pc}}
\def \Pt{\tilde{\P}}
\def \Ht{\tilde{H}}
\def \Hct{\tilde{\Hc}}
\def \Pb{\overline{\P}}
\def \Mb{\mathbf{M}}
\def \Pbb{\mathbf{P}}
\def \Pbc{\mathbf{P}^{\preceq}}
\def \Hcb{\overline{\Hc}}
\def \Hb{\overline{H}}
\def \Yb{\overline{Y}}

\def \lro{\longrightarrow}
\def \Th{\Theta}
\def \Zb{\overline{Z}}

\def\esup{\mathop{\rm ess \, sup}}

\title{Optimal Skorokhod embedding under finitely-many marginal constraints 
	\footnote{We are grateful to Jan Ob{\l}\'oj and three anonymous referees for helpful suggestions and comments.
	We gratefully acknowledge the financial support of the ERC 321111 Rofirm, the ANR Isotace, and the Chairs Financial Risks (Risk Foundation, sponsored by Soci\'et\'e G\'en\'erale) and Finance and Sustainable Development (IEF sponsored by EDF and CA).}}

\author{
	Gaoyue Guo\thanks{CMAP, Ecole Polytechnique, France. guo@cmap.polytechnique.fr}
	\and Xiaolu Tan\thanks{University of Paris-Dauphine, PSL Research University, CNRS, UMR [7534], CEREMADE.
	tan@ceremade.dauphine.fr}
	\and Nizar Touzi\thanks{CMAP, Ecole Polytechnique, France. nizar.touzi@polytechnique.edu}
}

\date{\today}

\begin{document}
\bibliographystyle{plain}

\maketitle

\abstract{The \textit{Skorokhod embedding problem} aims to represent a given probability measure on the real line as the distribution of Brownian motion stopped at a chosen stopping time. 
	In this paper, we consider an extension of the weak formulation of the \textit{optimal Skorokhod embedding problem} in Beiglb\"ock, Cox \& Huesmann \cite{BH} to the case of finitely-many marginal constraints
	\footnote{While producing the final version of this paper, we knew from Mathias Beiglb\"ock about the new development in \cite{BCHPP} that is also in the final producing version and extends the previous work \cite{BH} to the case of finitely-many marginal constraints.
	We emphasize that our approach is of completely different nature,
	and our results are established under more general conditions,
	although the dual formulations are slightly different in the two papers (see Section \ref{subsec:discussions} for more details).}.
	Using the classical convex duality approach together with the optimal stopping theory, 
	we establish some duality results under more general conditions than \cite{BH}.
	We also relate these results to the problem of martingale optimal transport under multiple marginal constraints. }

\vspace{2mm}


\vspace{2mm}





\section{Introduction}\label{sec:intro}

	Let $\mu$ be a probability measure on $\R$, with finite first moment and centered,
	the Skorokhod embedding problem (SEP) consists in finding a stopping time $\tau$ on a Brownian motion $W$ such that $W_{\tau} \sim \mu$ 
	and the stopped process $W_{\tau \wedge \cdot}:=\big(W_{\tau\wedge t}\big)_{t\ge 0}$ is uniformly integrable. We refer the readers to the survey paper \cite{Obloj} of Ob{\l}\'oj for a comprehensive account of the field.

	In this paper, we consider its extension to the case of multiple marginal constraints. Namely, let $\mu:=(\mu_1, \cdots, \mu_n)$ be a given family of centered probability measures such that the family is increasing in convex ordering, 
	i.e. for every convex function $\phi: \R \to \R$, one has
	\b*
		\int_{\R} \phi(x) \mu_k(dx)
		~\le~
		\int_{\R} \phi(x) \mu_{k+1}(dx)
		&\mbox{for all}~ k = 1, \cdots, n-1.
	\e*
	The extended SEP is to find an increasing family of stopping times $\tau:=(\tau_1, \cdots, \tau_n)$ 
	such that $W_{\tau_k} \sim \mu_k$ for all $k = 1, \cdots, n$
	and the stopped process $W_{\tau_n \wedge \cdot}$ is uniformly integrable.
	We study an associated optimization problem,   
	which consists in maximizing the expected value of some reward function among all such embeddings.

 One of the motivations to study this problem is its application in finance to compute the arbitrage-free model-independent price bounds of contingent claims consistent with the market prices of Vanilla options. Mathematically, the underlying asset is required to be a martingale according to the no-arbitrage condition  and the market calibration allows to recover the marginal laws of the underlying at certain maturities (see e.g. Breeden \& Litzenberger \cite{BL}).
Then by considering all martingales fitting the given marginal distributions, one can obtain the arbitrage-free price bounds. Based on the fact that every continuous martingale can be considered as a time-changed Brownian motion by Dambis-Dubins-Schwarz theorem, Hobson studied the model-free hedging of lookback options in his seminal paper \cite{Hobson1998} by means of the SEP. The main idea of his pioneering work is to exploit some solution of the SEP satisfying some optimality criteria, which yields the model-free hedging strategy and allows to solve together the model-free pricing and hedging problems. Since then, the optimal  SEP has received substantial attention from the mathematical finance community and various extensions were achieved in the literature, such as Cox \& Hobson \cite{CH}, Hobson \& Klimmek \cite{HK}, Cox, Hobson \& Ob{\l}\'oj \cite{CHO1}, Cox \& Ob{\l}\'oj \cite{CO1} and 
Davis, Ob{\l}\'oj \& Raval\cite{DOR},  Ob{\l}\'oj \& Spoida \cite{OblojSpoida}, etc. A thorough literature is provided in Hobson's survey paper \cite{Hobson2011}.

Beiglb\"ock, Cox and Huesmann generalized this heuristic idea and formulated the optimal SEP in \cite{BH}, which recovered many previous known results by a unifying formulation.
Namely, their main results are twofold. First, they establish the expected identity between the optimal SEP and the corresponding model-free superhedging problem. Second, they derive the characterization of the optimal embedding by a geometric pathwise property which allows to recover all previous known embeddings in the literature.

The problem of model-free hedging has also been approached by means of the martingale optimal transport, as initiated by Beiglb\"ock, Henry-Labord\`ere \& Penkner \cite{BHP} in the discrete-time case and Galichon, Henry-Labord\`ere \& Touzi \cite{GHT} in the continuous-time case. Further development enriches this literature, such as Beiglb\"ock \& Juillet \cite{BeiglbockJuillet}, , Henry-Labord\`ere \& Touzi \cite{HT}, Henry-Labord\`ere, Tan \& Touzi \cite{HTT}, etc. A remarkable contribution for the continuous-time martingale optimal transport is due to  Dolinsky \& Soner \cite{DS, DS2}. We also refer to Tan \& Touzi \cite{TanTouzi} for the optimal transport problem under more general controlled stochastic dynamics. 

Our objective of this paper is to revisit the duality result of \cite{BH} and to extend the duality under more general conditions and to the case of multiple marginal constraints. Our approach uses tools from a completely different nature. First, by following the convex duality approach, we convert the optimal SEP into an infimum of classical optimal stopping problems. Next, we use the standard dynamic programming approach to relate such optimal stopping problems to model-free superhedging problems. 
We observe that the derived duality allows to reproduce the geometric characterization of the optimal embedding introduced in \cite{BH}, see e.g. \cite{GTT2}. 
Finally, we show that our result induces the duality for a class of martingale optimal transport problems in the space of continuous paths.

The paper is organized as follows. In Section \ref{sec:osep}, we formulate our optimal SEP under finitely-many marginal constraints and provide two duality results. In Section \ref{sec:MT}, the duality of optimal SEP together with time-change arguments gives the duality for the martingale optimal transport problem under multi-marginal constraints. We finally provide the related proofs in Section \ref{sec:proofs}. 

	\vspace{2mm}

	\noindent {\bf Notations.}
	\rmi Let $\Om := C(\R_+, \R)$ be the space of all continuous paths $\om$ on $\R_+$ such that $\om_0 = 0$, 
	$B$ be the canonical process, $\P_0$ be the Wiener measure,
	$\F := (\Fc_t)_{t \ge 0}$ be the canonical filtration generated by $B$,
	and $\F^a := (\Fc^a_t)_{t \ge 0}$ be the augmented filtration under $\P_0$.
	
	\vspace{1mm}

	\noindent \rmii Define for some fixed integer $n \ge 1$ the enlarged canonical space by $\Omb := \Om \x \Th$ (see El Karoui \& Tan \cite{EKT1, EKT2}), where $\Th :=\big\{ (\theta_1, \cdots, \theta_n) \in \R_+^n: \theta_1 \le \cdots \le \theta_n \big\}$. All the elements of $\Omb$ are denoted by $\omb:= (\om,\theta)$ 
	with $\theta :=(\theta_1,\cdots, \theta_n)$. Denote further by $(B, T)$ (with $T:=(T_1,\cdots,T_n )$) the canonical element on $\Omb$, i.e. $B_t(\omb):=\om_t$ and $T(\omb):=\theta$ for every $\omb=(\om,\theta)\in\Omb$.
	The enlarged canonical filtration is denoted by $\Fbb:= (\Fcb_t)_{t \ge 0}$,
	where $\Fcb_t$ is generated by $(B_s)_{0\le s \le t}$ and all the sets $\{ T_k \le s \}$ for all $s\in [0,t]$ and $k=1,\cdots, n$. 
	In particular, all random variables $T_1, \cdots, T_n$ are $\Fbb-$stopping times.

	\vspace{1mm}

	\noindent \rmiii We endow
	$\Om$ with the compact convergence topology,
	and $\Theta$ with the classical Euclidean topology,
	then $\Om$ and $\Omb$ are both Polish spaces (separable, complete metrizable space).
	In particular, $\Fcb_{\infty} := \bigvee_{t \ge 0} \Fcb_t$ is the Borel $\sigma-$field of the Polish space $\Omb$
	(see Lemma \ref{lemm:Fcb}).

	\vspace{1mm}

	\noindent \rmiv Denote by $\Cc_1:=\Cc_1(\R)$ the space of all continuous functions on $\R$ with linear growth.

	\vspace{1mm}

\noindent \rmv Throughout the paper UI, a.s. and q.s. are respectively the abbreviations of uniformly integrable, almost surely and quasi-surely. Moreover, given a set of probability measures $\Nc$ (e.g. $\Nc=\Pcb$ and $\Nc=\Mc$ in the following) on some measurable space, we write $\Nc-$q.s. to represent that some property holds under every probability of $\Nc$.

\section{An optimal Skorokhod embedding problem and the dualities}
\label{sec:osep}

	In this section, we formulate an optimal Skorokhod embedding problem (SEP) under finitely-many marginal constraints, as well as its dual problems. We then provide two duality results.

\subsection{An optimal Skorokhod embedding problem}

	Throughout the paper, $\mu := (\mu_1, \cdots, \mu_n)$ is a vector of $n$ probability measures on $\R$ and we denote, for any integrable function $\phi: \R\to\R$,
\b*
	\mu_k(\phi)&:=&\int_{\R} \phi(x)\mu_k(dx) ~~\mbox{for all } k=1,\cdots, n.
\e*
The vector $\mu$ is said to be a peacock if each probability $\mu_k$ has finite first moment, i.e. $\mu_k(|x|)<+\infty$, and $\mu$ is increasing in convex ordering, i.e, $k\mapsto \mu_{k}(\phi)$ is non-decreasing for every convex function $\phi$. A peacock $\mu$ is called centered if $\mu_k(x)=0$ for all $k=1, \cdots, n$. 
Denote by $\Pbb^{\preceq}$ the collection of all centered peacocks. 

\paragraph{Optimal SEP}
	As in Beiglb\"ock, Cox \& Huesmann \cite{BH}, we shall consider the problem in a weak setting, 
	i.e. the stopping times may be identified by probability measures on the enlarged space $\Omb$.
	Recall that the elements of $\Omb$ are denoted by $\omb:=\big(\om,\theta=(\theta_1,\cdots,\theta_n)\big)$ and the canonical element is denoted by $\big(B, T=(T_1,\cdots, T_n)\big)$, and in particular $T_1, \cdots, T_n$ are all $\Fbb-$stopping times. Let $\Pcb(\Omb)$ be the space of all probability measures on $\Omb$,
	and define
	\be \label{eq:def_Pcb}
		\Pcb
		~:=~
		\Big\{
			\Pb \in \Pcb(\Omb):
			B ~\mbox{is an}~ \Fbb-\mbox{Brownian motion and }
			B_{T_n \wedge \cdot}
			~\mbox{is UI under } \Pb\Big\}.
	\ee
	Set for any given family of probability measures $\mu = (\mu_1, \cdots, \mu_n)$
	\be\label{def:embedding}
		\Pcb (\mu)
		&:=&
		\Big\{
			\Pb \in \Pcb:
			B_{T_k} \stackrel{\Pb}{\sim}\mu_k \mbox{ for all } k = 1, \cdots, n
		\Big\}.
	\ee
	As a consequence of Kellerer's theorem in \cite{Kellerer}, $\Pcb(\mu)$ is nonempty if and only if 
	$\mu\in\Pbb^{\preceq}$.

\vspace{2mm}

Let $\Phi : \Omb \to \R$ be a Borel measurable function, then $\Phi$ is called non-anticipative if 
$\Phi (\om, \theta)=\Phi \big(\om_{\theta_n \wedge \cdot}, \theta\big)$ for every $(\om,\theta)\in \Omb$.
Define the optimal SEP for a non-anticipative function $\Phi$ by
	\be \label{eq:P}
		P(\mu)
		&:=&
		\sup_{\Pb \in \Pcb(\mu)}
		~\E^{\Pb} 
		\big[\Phi(B,T)\big],
	\ee
where the expectation of a random variable $\xi$ is defined by $\E^{\Pb}[ \xi] = \E^{\Pb}[ \xi^+] - \E^{\Pb}[ \xi^-]$ with the convention $\infty - \infty = - \infty$.
The problem is well-posed if there exists at least a $\Pb\in\Pcb(\mu)$ such that $\E^{\Pb}[|\Phi(B,T)|]<+\infty$. We emphasize that $\Phi$ is assumed to be non-anticipative throughout the paper.

\begin{Remark}
	\rmi	A $\mu-$embedding is a collection
	\b*
		\alpha 
		&=& 
		\big( \Om^{\alpha}, \Fc^{\alpha}, \P^{\alpha}, \F^{\alpha} = (\Fc^{\alpha}_t)_{t \ge 0}, W^{\alpha}, \tau^{\alpha}=(\tau^{\alpha}_1, \cdots, \tau^{\alpha}_n) \big),
	\e*
	where $W^{\alpha}$ is an $\F^{\alpha}-$Brownian motion,
		$\tau^{\alpha}_1, \cdots, \tau^{\alpha}_n$ are increasing  $\F^{\alpha}-$stopping times such that 			$ W^{\alpha}_{\tau^{\alpha}_n\wedge\cdot}$ is uniformly integrable,
		and $W^{\alpha}_{\tau_k^{\alpha}}\sim\mu_k$ for all $k=1,\cdots n$.
		We observe that for every centered peacock $\mu$, every $\mu-$embedding $\alpha$ induces a probability measure $\Pb:=\P^{\alpha}\circ (W^{\alpha}, \tau^{\alpha})^{-1}\in\Pcb(\mu)$.
		Conversely,  every probability measure $\Pb\in\Pcb(\mu)$ together with the canonical space $(\Omb, \Fcb_{\infty})$, canonical filtration $\Fbb$, and canonical element $(B, T)$ is a $\mu-$embedding.
		Then denoting by $\Ac(\mu)$ the collection of all $\mu-$embeddings,
		the optimal SEP \eqref{eq:P} is equivalent to
		$$
			\sup_{\alpha \in \Ac(\mu)} \E^{\P^{\alpha}}\big[ \Phi(W^{\alpha}, \tau^{\alpha}) \big].
		$$

\noindent \rmii The problem \eqref{eq:P} can be considered as a weak formulation of the optimal SEP.
	A strong formulation consists in considering all stopping times w.r.t. the Brownian filtration,
	and it may not be equivalent to the weak formulation
	(especially when $\mu$ has an atom at $0$, see Example \ref{exam:non_equiv}).
	Although most of the well known embeddings are ``strong'' stopping times,
	some optimal embeddings are constructed in ``weak'' sense, such as that in Hobson \& Pedersen \cite{HP}.
	We also notice that it should be natural to consider the weak formulation to obtain the existence of the optimizer in general cases, since the space of all ``weak'' embeddings is compact under the weak convergence topology as shown below.
\end{Remark}

\subsection{The duality results}

	We introduce two dual problems. Recall that $\P_0$ is the Wiener measure on $\Om=C(\R_+, \R)$ under which the canonical process $B$ is a standard Brownian motion, 
	$\F = (\Fc_t)_{t \ge 0}$ is the canonical filtration 
	and $\F^a = (\Fc^a_t)_{t \ge 0}$ is the  $\P_0-$augmented filtration.
	Denote by $\Tc^a$ the collection of all increasing families of $\F^a-$stopping times $\tau = (\tau_1, \cdots, \tau_n)$ such that the process $B_{\tau_n \wedge \cdot}$ is uniformly integrable.
	Define also the class of functions
	\be\label{eq:Lambda}
		\Lambda
		&:=&
		\Cc_1^n
		~=~
		\Big\{
			\lambda := (\lambda_1, \cdots, \lambda_n)
			~: \lambda_k \in \Cc_1 \mbox{ for all } k=1, \cdots n
		\Big\}.
	\ee 
For $\mu=(\mu_1,\cdots,\mu_n)$, $\lambda=(\lambda_1,\cdots, \lambda_n)$ and $\big(\om,\theta=(\theta_1,\cdots,\theta_n)\big)\in\Omb$, we denote
	\b*
	\mu(\lambda)
		~:=~
		\sum_{k =1}^n \mu_k(\lambda_k) 	&\mbox{and}& 
\lambda(\om_{\theta})~:=~\sum_{k =1}^n \lambda_k(\om_{\theta_k})
		\mbox{ with } \om_{\theta}~:=~(\om_{\theta_1},\cdots,\om_{\theta_n}).
	\e*
	Then the first dual problem for the optimal SEP \eqref{eq:P} is given by
	\be \label{eq:D1}
		D_0(\mu)
		&:=&
		\inf_{\lambda \in \Lambda} \Big\{
			\sup_{\tau \in \Tc^a} \E^{\P_0}
			\big[
				\Phi (B,\tau) 
				-\lambda(B_{\tau}) 
			\big] 
			+ \mu(\lambda)
		\Big\}.
	\ee
	As for the second dual problem, we return to the enlarged space $\Omb$. Given $\Pb\in\Pcb$, an $\Fbb-$optional process $M = (M_t)_{t \ge 0}$ is called a strong $\Pb-$supermartingale if
\b*
		\E^{\Pb} \big[ M_{\tau_2}\big |\Fcb_{\tau_1} \big]
		&\le&
		M_{\tau_1},
		~\Pb-\mbox{a.s.}
	\e*	
	for all $\Fbb-$stopping times $\tau_1 \le \tau_2$. Let $\L^2_{loc}$ be the space of all $\Fbb-$progressively measurable processes $\Hb =(\Hb_t)_{t \ge 0}$ such that
	\b*
		\int_0^{t} \Hb_s^2 ds~<~ +\infty \mbox{ for every } t\ge 0, ~\Pcb-\mbox{q.s.}.
	\e*
	For $\Hb\in\L^2_{loc}$, the stochastic integral $(\Hb\cdot B):=\int_0^{\cdot}\Hb_sdB_s$ is well defined $\Pb-$a.s. for all $\Pb\in\Pcb$. We introduce a subset of processes:
	\b*
		\Hcb
		&:=&
		\Big \{ \Hb \in \L^2_{loc}~:  
			(\Hb\cdot B)
			~\mbox{is a}~ \P-\mbox{strong supermartingale for all}~ \Pb \in \Pcb
		\Big\}.
	\e*
	Denote further
	\b*
		\Dc &:=&
		\Big\{ 
			(\lambda, \Hb) \in \Lambda \x \Hcb
			:
			\lambda(B_T) +(\Hb\cdot B)_{T_n}
			\ge \Phi (B,T),~ \Pcb-\mbox{q.s.}
		\Big\},
	\e*
	and the second dual problem is given by
	\be \label{eq:D2}
		D(\mu)
		&:=&
		\inf_{(\lambda, \Hb) \in \Dc}\mu(\lambda).
	\ee
Loosely speaking, the two dual problems dualize respectively different constraints of the primal problem \eqref{eq:P}. By penalizing the marginal constraints, we obtain the first dual problem $D_0(\mu)$ of \eqref{eq:D1}, where a multi-period optimal stopping problem appears for every fixed $\lambda \in \Lambda$. Then the second dual problem $D(\mu)$ of \eqref{eq:D2} follows by the resolution of the optimal stopping problem via the Snell envelope approach and the Doob-Meyer decomposition.

	Our main duality results require the following conditions.

	\begin{Assumption} \label{Hyp:phi}
		The reward function $\Phi: \Omb \to \R$ is Borel measurable, non-anticipative,
		bounded from above, and
		$\theta \mapsto \Phi(\om_{\theta_n \wedge \cdot}, \theta)$ is upper-semicontinuous for 
		$\P_0-$a.e. $\om \in \Om$.
	\end{Assumption}

	\begin{Assumption} \label{Hyp:phi2}
		One of the following conditions holds true.
		
		\noindent \rmi $n=1$.
		
		\noindent \rmii $n \ge 2$ and the map $\omb \mapsto \Phi (\omb)$ is upper-semicontinuous.

		\noindent \rmiii $n \ge 2$ and the reward function $\Phi$ admits the representation 
		$$
			\Phi (\omb) 
			~~=~~
			\sum_{k=1}^n \Phi_k (\om, \theta_1, \cdots, \theta_k),
		$$
		where for each $k = 1, \cdots, n$, $\Phi_k : \Om \x \R_+^k \to \R$ satisfies that
		$\Phi_k(\om, \theta_1, \cdots, \theta_k) = \Phi_k(\om_{\theta_k \wedge \cdot}, \theta_1, \cdots, \theta_k)$,
		and $(\theta_1, \cdots, \theta_{k-1}) \mapsto \Phi_k (\om_{\theta_k \wedge \cdot}, \theta_1, \cdots, \theta_k)$
		is uniformly continuous for $0 \le \theta_1 \le \cdots \le \theta_{k-1} \le \theta_k$, 
		uniformly in $\theta_k$.

	\end{Assumption}

	\begin{Theorem} \label{theo:main}
		\rmi Under Assumption \ref{Hyp:phi}, there is some $\Pb^* \in \Pcb(\mu)$ such that 
		\b*
		\E^{\Pb^*}[ \Phi ] ~=~ P(\mu) ~=~ D_0 (\mu).
		\e*
		\noindent \rmii Suppose in addition that Assumption \ref{Hyp:phi2} holds true, then
		\b*
			P(\mu) ~~=~~ D_0(\mu) ~~=~~ D(\mu).
		\e*
	\end{Theorem}

\paragraph{The case of a separable reward function}

	When $\Phi$ is of the form introduced in Assumption \ref{Hyp:phi2} $\mathrm{(iii)}$,
	we can consider a stronger dual formulation.
	Denote by $\Hc$ the collection of all $\F-$predictable processes $H^0: \R^+ \x \Om \to \R$ 
	such that the stochastic integral $(H^0 \cdot B)_t := \int_0^t H^0_s dB_s$ is a martingale under $\P_0$,
	and $(H^0 \cdot B)_t \ge - C(1 + |B_t|)$ for some constant $C > 0$.
	
	In the filtered space $(\Om, \Fc, \P_0, \F)$,
	we say a process $X$ is of class (DL) if for each $t \ge 0$, 
	the family $\{ X_{\tau} ~: \tau \le t ~\mbox{is a stopping time}\}$ is uniformly integrable;
	we say an $\F-$optional process $X$ of class (DL) is an $\F-$ supermartingale if for all bounded stopping times $\sigma \le \tau$, one has $X_{\sigma} \ge \E^{\P_0}[ X_{\tau} | \Fc_{\sigma}]$.
	Denote further by $\Sc$ the set of all $\F-$supermartingales in $(\Om, \Fc, \P_0)$ such that
	$|S_t| \le C(1 + |B_t|)$ for some constant $C > 0$.

	Define then
	\b*
		\Dc'
		\!\!\!&:=&\!\!\!
		\Big\{
		(\lambda, H^1, \cdots, H^n) \in \Lambda \x (\Hc)^n
		~:
		\sum_{k=1}^n \big(  \lambda_k(\om_{\theta_k}) +  \int_{\theta_{k-1}}^{\theta_k} H^k_s d B_s\big)
		\ge \Phi \big(\om, \theta \big) \\
		&&
		~~~~~~~~~~~~~~~~~~~~~~~~~~~~~~~~
		\mbox{for all}~ 0 \le \theta_1 \le \cdots \le \theta_n,
		~\mbox{and}~\P_0 -\mbox{a.e.}~ \om \in \Om
		\Big\},
	\e*
	and
	\b*
		\Dc''
		\!\!\!&:=&\!\!\!
		\Big\{
		(\lambda, S^1, \cdots, S^n) \in \Lambda \x (\Sc)^n
		~:
		\sum_{k=1}^n \big(  \lambda_k(\om_{\theta_k}) + S^k_{\theta_k}  - S^k_{\theta_{k-1}} \big)
		\ge \Phi \big(\om, \theta \big) \\
		&&
		~~~~~~~~~~~~~~~~~~~~~~~~~~~~~~~~
		\mbox{for all}~ 0 \le \theta_1 \le \cdots \le \theta_n,
		~\mbox{and}~\P_0 -\mbox{a.e.}~ \om \in \Om
		\Big\}.
	\e*

	\begin{Proposition} \label{prop:main_p}
		Suppose that Assumption \ref{Hyp:phi} and Assumption \ref{Hyp:phi2} \rmiii hold true.
		Suppose in addition that $\Phi_k (\om, \theta_1, \cdots, \theta_k)$ depends only on $(\om, \theta_k)$.
		Then
		\be \label{eq:D_prime}
			P(\mu)
			~~=~~
			D'(\mu)
			~~:=~~
			\inf_{(\lambda, H) \in \Dc'} \mu(\lambda)
			~=~
			D''(\mu)
			~~:=~~
			\inf_{(\lambda, S) \in \Dc''} \mu(\lambda).
		\ee
	\end{Proposition}

\subsection{More discussions and examples}
\label{subsec:discussions}

	\begin{Remark} \label{rem:compar_BCH}
		The above dual formulation \eqref{eq:D2} has been initially provided and proved by \cite{BH} in the one marginal case  ($n=1$),
		under the condition that $(\om, \theta) \mapsto \Phi(\om_{\theta \wedge \cdot}, \theta)$ is bounded from above and upper semicontinuous.
		
		When $n=1$, our duality results hold under more general conditions: $\Phi$ is non-anticipative, bounded from above, and $\theta \mapsto \Phi(\om_{\theta \wedge \cdot}, \theta)$ is u.s.c. for $\P_0-$a.e. $\om \in \Om$.
		In particular, this allows to include the case where $\Phi$ is a function of the local time of the stopped Brownian motion, since the local time of the Brownian motion is continuous in $\theta$ but has no regularity in $\om$.
		As an important example, the optimal embedding w.r.t. a convex function of the local time is provided by Vallois's embedding, see e.g. \cite{CHO1, CGH}.

		Nevertheless, we use a quasi-sure formulation in our dual problem \eqref{eq:D2}, and a pathwise formulation in \eqref{eq:D_prime}.
		The dual problem in \cite{BH, BCHPP} uses a pathwise formulation.
		Moreover, instead of the stochastic integral $(H \cdot B)$ in our case,
		they use martingales which are continuous in $(t,\om)$ in the dual formulation.

		For the multiple marginal case ($n\ge 2$), 
		when $\Phi$ has no regularity in $\om$,
		we need a uniform continuity condition in time variables $(\theta_1, \cdots, \theta_{n-1})$ but not in $\theta_n$.
		The uniform continuity condition is a technical condition to aggregate a family of supermartingales appearing in the classical optimal stopping problem.
		We can next approximate an u.s.c. function by a sequence of Lipschitz functions.
		However, to keep the non-anticipative property of $\Phi$, we need to assume that $\Phi$ is u.s.c. w.r.t. both variables $(\om, \theta)$ in Assumption \ref{Hyp:phi2} $\mathrm{(ii)}$ (see the proof in Section \ref{subsubsec:reduction}).
		This is also the main reason for the regularity conditions in $\om$ in K\"allblad, Tan \& Touzi \cite{KTZ},
		where the duality result is extended to the infinitely-many marginals case.
	\end{Remark}

	\begin{Remark} 
		 A characterization of the optimizers $\Pb^*$ has been provided in \cite{BH},
		called monotonicity principle.
		An alternative proof of this result is given in our accompanying paper \cite{GTT2}.

		\vspace{1mm}

		\noindent \rmii For the general martingale optimal transport problem, 
		the dual optimizer $\lambda^*$ may not exist, 
		as shown in Beiglb\"ock, Henry-Labord\`ere \& Penkner \cite{BHP}.
		More recently, by relaxing the dual formulation of the one dimensional discrete-time martingale transport, the existence of the dual optimizer in ``weak'' sense is  obtained by Beiglb\"ock, Nutz \& Touzi \cite{BNT}.
		The existence of dual optimizer in our context is still an open question.

		\noindent \rmiii Nevertheless, when the function $\Phi$ has a particular form, we do have the dual optimizer $\lambda^*$, and the optimizers $\Pb^*$ and $\lambda^*$ can be explicitly constructed.
		For example, Hobson \cite{Hobson1998} provided the construction when $n=1$ and $\Phi$ is an increasing function of the running maximum,
		Hobson \& Klimmek \cite{HK} studied the case for the forward starting straddle,
		Cox, Hobson \& Ob{\l}\'oj \cite{CHO1} considered functions on local time,
		see also Brown, Hobson \& Rogers \cite{BHR1, BHR2},
Cox \& Ob{\l}\'oj \cite{CO1}, Davis, Ob{\l}\'oj \& Raval\cite{DOR}, etc. among many others, for more concrete cases.
		We also refer to Hobson \cite{Hobson2011} for a detailed review on these constructions.
	\end{Remark}

	\begin{Remark} 		
		Based on the first dual problem $D_0(\mu)$ in \eqref{eq:D1}, 
		a numerical algorithm has been obtained in Bonnans and Tan \cite{BT} for the above optimal SEP.
	\end{Remark}

	\begin{Remark}
		To prove the equality $D_0(\mu) = D(\mu)$, we study a multiple optimal stopping problem using a backward iteration approach, since the stopping times $T_1, \cdots, T_n$ are assumed to be ordered.
		The order condition $T_1 \le \cdots \le T_n$ is natural as motivated by its applications in finance (see Section \ref{sec:MT}) and technically necessary in our arguments.
		Without the order condition, one can always formulate an optimal SEP, but the corresponding dual problem seems not clear.
	\end{Remark}

	\begin{Example} \label{exam:Phi}
		\rmi Let $\phi: \R_+ \x (\R^3)^n$ be a continuous function, bounded from above,
		denote $\overline \om_t := \sup_{0 \le s \le t} \om_s$ and $\underline \om_t := \inf_{0 \le s \le t} \om_s$.
		Since $\om \mapsto (\om_t, \overline \om_t, \underline \om_t)$ is continuous,
		the reward function $\Phi$ defined by
		$$
			\Phi(\om, \theta) 
			~:=~
			\phi \big(\theta_i, \om_{\theta_i}, \overline \om_{\theta_i}, \underline \om_{\theta_i}, i=1, \cdots, n \big)
		$$
		satisfies clearly Assumptions \ref{Hyp:phi} and \ref{Hyp:phi2} $\mathrm{(ii)}$.

		\vspace{1mm}

		\noindent \rmii Let $L: \Om \x \R_+ \to \R$ be the local time of the Brownian motion.
		We can choose $L$ to be $\F-$predictable since any $\F^a-$predictable process is indistinguishable to  an $\F-$predictable process.
		Then $t \mapsto L_t(\om)$ is continuous and increasing for $\P_0-$a.e. $\om \in \Om$.
		Let $\phi: \R_+ \to \R$ be a continuous function, bounded from above,
		then $\Phi(\omb) := \phi(L_{\theta_n}(\om))$ satisfies Assumptions \ref{Hyp:phi} and \ref{Hyp:phi2} $\mathrm{(iii)}$.
	\end{Example}

	\begin{Example}[Nonequivalence between the strong and weak formulation]
	\label{exam:non_equiv} 
		When $\mu$ has an atom at $0$, one can easily show the nonequivalence between the strong and weak formulation.
		Let $n=1$, $\mu := \frac13 \delta_{\{0\}} +  \frac13 \delta_{\{1\}} +  \frac13 \delta_{\{-1\}}$ and $\Phi(\om, \theta) := \1_{\{0\}}(\theta)$.
		Define $\tau_0 := \inf\{ t ~:|B_t| \ge 1\}$, and $\Pb_0 := \frac13 \P_0 \circ (B,0)^{-1} + \frac23 \P_0 \circ (B, \tau_0)^{-1}$,
		then $\Pb_0 \in \Pc(\mu)$ and $\E^{\Pb} \big[ \Phi(B, T) \big] = \frac13$.
		Further, let $\tau \in \Tc^a$ such that $B_{\tau} \sim \mu$ under the Wiener measure $\P_0$,
		then $\P_0[\tau > 0] > 0$.
		Since the
		augmented Brownian filtration satisfies Blumenthal's zero-one law, then $\P_0[\tau > 0] = 1$.
		It follows that
		\b*
			\sup_{\tau \in \Tc^a, B_{\tau} \sim \mu} \E^{\P_0} \big[ \Phi(B, \tau) \big]
			~=~
			0
			~<~
			\frac13
			~\le~
			\sup_{\Pb \in \Pcb(\mu)} \E^{\Pb} \big[ \Phi(B, T) \big].
		\e*
		
	\end{Example}

	We finally provide an example where the duality fails when $\Phi(\om, \theta)$ has no regularity in $\theta$.

	\begin{Example}
		Let $n=1$,  $\Phi(\om, \theta) := \1_{\Q}(\theta)$, where $\Q$ denotes the set of all rational numbers, and $\mu := \frac12 \delta_{\{1\}} + \frac12 \delta_{\{-1\}}$.
		We first notice that $\Pcb(\mu)$ has only one element, which is the probability measure induced 
		by $(B, \tau_0)$, where $B$ is a standard Brownian motion and $\tau_0 := \inf\{ t ~:|B_t| \ge 1\}$.
		Indeed, for any $\Pb \in \Pcb(\mu)$, one has $\E^{\Pb}[T] = \E^{\Pb}[B_T^2] = \E^{\Pb}[\tau_0(B)]$
		and $T \ge \tau_0(B)$, $\Pb-$a.s.
		Moreover, since the hitting time $\tau_0$ is a random variable of continuous distribution on $\R_+$,
		then
		\b*
			P(\mu)~ =~ \sup_{\Pb \in \Pcb(\mu)} \E^{\Pb} \big[ \Phi(B, T) \big]
			~=~ \E^{\Pb_0} \big[ \1_{\Q} (\tau_0(B)) \big]
			~=~ 0.
		\e*
		 As for the dual problem, we notice that $\lambda \in \Lambda$ is a continuous function,
		and one can approximate a stopping time by stopping times taking value in $\Q$, then
		\b*
			\sup_{\tau \in \Tc^a} \E^{\P_0} \big[ \1_{\Q}(\tau) - \lambda (B_{\tau}) ]
			~=~
			\sup_{\tau \in \Tc^a} \E^{\P_0} \big[ 1 - \lambda (B_\tau) ],
			~~\mbox{for all}~ \lambda \in \Lambda.	
		\e*
		Then by its definition in \eqref{eq:D1}, 
		$D_0(\mu) = \inf_{\lambda \in \Lambda} \big\{\mu(\lambda) +\sup_{\tau \in \Tc^a} \E^{\P_0} \big[ 1 - \lambda (B_\tau) ] \big\} = 1$.
		Similarly, we can easily deduce that for $\Phi(\omb) = \1_{\Q}(\theta)$,
		$D(\mu) = 1$ and  it follows that
		\b*
			P(\mu) ~=~ 0 ~\neq 1 ~=~ D_0(\mu) ~=~ D(\mu),
		\e*
		in the above context.
	\end{Example}

\section{Application to a class of martingale transport problems}
\label{sec:MT}

	In this section,
	we use the previous duality results of the optimal SEP 
	to study a continuous time martingale transport problem under multiple marginal constraints.
	As an application in finance to study the robust superhedging problem,
	the multi-marginal case is very natural.
	Namely, when the Vanilla options are available for trading for several maturities, 
	thus inducing
	the marginal distributions of the underlying asset at several times,
	we can formulate the robust superhedging problem as a martingale transport problem under multiple  marginal constraints.

\subsection{Robust superhedging and martingale transport}

Define the canonical process $X:=(X_t)_{0\le t\le 1}$ by $X_t=B_{1\wedge t}$ for all $t\in [0,1]$ and its natural filtration $\Ft := (\Fct_t)_{0 \le t \le 1}$. Denote further by $\Mc$ the collection of all martingale measures $\Pt$, i.e. the probability measures under which $X$ is a martingale. Let $I := (0<t_1<\cdots<t_n=1)$ be a set of time instants and define the set of martingale transport plans for $\mu\in\Pbb^{\preceq}$
	\b*
		\Mc(\mu) 
		&:=&
		\Big\{ \Pt \in \Mc ~: X_{t_k} \stackrel{\Pt}{\sim} \mu_k \mbox{ for all } k = 1, \cdots, n\Big\}.
	\e*
	By Karandikar \cite{Karandikar}, there is a non-decreasing $\Ft-$progressive process $\Xcr$
	taking value in $[0,\infty]$,
	such that $\Xcr$ coincides with the quadratic variation of $X$, $\Pt-$a.s.
	for every martingale measure $\Pt \in \Mc$.
	Denote $\Xcr^{-1}_t := \inf \big\{ s\ge 0 ~: \Xcr_s > t\big\} \wedge 1$ and
	\be \label{eq:TimeChangeBM}
		W_t
		~:=~
		X_{\Xcr^{-1}_t} \1_{\{t < \Xcr_1\}} +  \big( X_1 + \widehat W_{t- \Xcr_1}\big) \1_{\{t \ge \Xcr_1\}},
	\ee
	where $\widehat W$ is an independent Brownian motion
	\footnote{In general case, one needs to enlarge the space to obtain an independent Brownian motion $\widehat W$. 
	However, in the following we will always consider a non-anticipative functional $\Phi(W_{\Xcr_{t_n}\wedge\cdot}, \Xcr_{t_1}, \cdots, \Xcr_{t_n})$ for $0 \le t_1 \le \cdots \le t_n=1$, which does not really depend on $\widehat W$.
	}.
	Then it follows from the Dambis-Dubins-Schwarz theorem (see e.g. Revuz \& Yor \cite[Theorem 1.7, Chapter V]{RevuzYor}) that the process $W$ is a Brownian motion.
	We denote also $W(X):=(X_{\Xcr^{-1}_t})_{0\le t\le \Xcr_1}$,
	which depends only on $X$.
	For a measurable function $\xi : \Om \to \R$, the martingale transport problem under multiple marginal constraints is defined by
	\be \label{eq:Pt}
		\tilde P (\mu)
		&:=&
		\sup_{\Pt \in \Mc(\mu)}
		\E^{\Pt} \big[ \xi (X) \big].
	\ee
Denote by $\Hct$ the collection of all $\Ft-$progressive processes $\Ht:=(\Ht_t)_{0\le t\le 1}$ such that
	\b*
		\int_0^1 \!\! \Ht_s^2 d \langle X \rangle_s <+\infty,~ \Mc-\mbox{q.s. and}~
		(\Ht \cdot X)
		~\mbox{is}~ \Pt-\mbox{supermartingale for all}~\Pt \in \Mc.
	\e*
	Then the two dual problems are given by	
	\be \label{eq:Dt}
		\tilde D_0(\mu)
		:=
		\inf_{\lambda \in \Lambda} \Big\{
			\sup_{\Pt \in \Mc} \E^{\Pt} \big[ 
				\xi \big(X\big) - \lambda(X_I)
			\big]
			+ \mu(\lambda)
		\Big\}
		~\mbox{and}~
	\tilde D(\mu)
		:=
		\inf_{(\lambda, \Ht) \in \tilde \Dc}\mu(\lambda), 
	\ee
	where
	\b*
		\lambda(X_I)&:=&\sum_{i=1}^n\lambda_i(X_{t_i}) \mbox{ with } X_I~:=~(X_{t_1},\cdots, X_{t_n})
	\e*
	and
	\b*
		\tilde \Dc
		&:=&
		\Big\{ (\lambda, \Ht) \in \Lambda \x \Hct: 
			\lambda(X_I) + (\Ht\cdot X)_1 \ge \xi \big(X\big),
			~\Mc-\mbox{q.s.}
		\Big\}.
	\e*
	It is easy to check that the weak dualities hold:
	\be \label{eq:weak_dual2}
		\tilde P (\mu)
		~~\le~~
		\tilde D_0(\mu)
		~~\le~~
		\tilde D(\mu).
	\ee

\subsection{Duality and financial interpretations}

	Using the duality results of the optimal SEP in Theorem \ref{theo:main}, we can establish the duality for the above martingale transport problem.
	
	\begin{Theorem} \label{theo:main2}
		Assume that the reward function $\xi$ admits the representation
		\b*
			\xi ( X)
			&=&
			\Phi \big(W, \Xcr_{t_1}, \cdots, \Xcr_{t_n}\big) \mbox{ with } W=W(X)
		\e*
		for some $\Phi: \Omb \to \R$ satisfying Assumptions \ref{Hyp:phi} and \ref{Hyp:phi2}. Then
		\b*
			\tilde P(\mu) 
			~~=~~ \tilde D_0(\mu)
			~~=~~ \tilde D(\mu).
		\e*
	\end{Theorem}

\paragraph{Financial interpretations}

	\begin{Example}
	Let $\phi: \big(\R_+ \x \R^3 \big)^n \to \R$ be a continuous function, bounded from above
	and $\xi$ be defined by
	\be \label{eq:xi}
		\xi(X)
		&=&
		 \phi \big( \Xcr_{t_i}, X_{t_i}, \overline{X}_{t_i}, \underline{X}_{t_i}, i = 1, \cdots n \big),
	\ee
	where $\overline{X}_t := \sup_{0 \le s \le t} X_s$ and $\underline{X}_t := \inf_{0 \le s \le t}X_s$.
	Then with 
		$$
			\Phi(\om, \theta) 
			~:=~
			\phi \big(\theta_i, \om_{\theta_i}, \overline \om_{\theta_i}, \underline \om_{\theta_i}, i=1, \cdots, n \big),
		$$		
	where $\overline \om_t := \sup_{0 \le s \le t} \om_s$ and $\underline \om_t := \inf_{0 \le s \le t} \om_s$,
	it is clear that $\xi$ satisfies the conditions in Theorem \ref{theo:main2}
	(see also Example \ref{exam:Phi}).
	The form \eqref{eq:xi} covers a big class of payoff functions of lookback option, barrier options, variance options, etc.
	\end{Example}

	The duality results in Theorem \ref{theo:main2} relates a problem of the arbitrage-free price bound with the 
	minimum robust superhedging problem.
	A martingale measure $\Pt \in \Mc$ can be considered as a market model,
	and the expectation of $\xi$ under a martingale measure provides an arbitrage-free price of option $\xi$.
	Then a probability measure $\Pt \in \Mc(\mu)$ can be considered 
	as a martingale model calibrated to the market information,
	since one can recover the marginal distribution $\mu$ of the underlying,
	when the Vanilla options at certain maturities are rich enough on the market (see e.g. \cite{BL}).
	Thus the primal problem \eqref{eq:Pt} provides an arbitrage-free price bound.

	As for the dual problem \eqref{eq:Dt}, 
	$\lambda$ and $\tilde H$ defines a semi-static strategy which superreplicates the payoff $\xi$ almost-surely under all possible martingale models.
	Then $D(\mu)$ provides the minimal robust superhedging cost of the exotic option $\xi$,
	using a class of possible static and dynamic strategies.
	Here robustness refers to the fact that the underlying probability measure is not fixed a priori, so that the superhedging requirement is imposed under all possible models $\Pt \in \Mc$.

	In Dolinsky \& Soner \cite{DS}, the duality is established (in a stronger sense) 
		for the case $n =1$, 
		for a general payoff function $\xi$
		which is Lipschtiz with respect to the uniform metric.
		In our Theorem \ref{theo:main2}, the reward function $\xi$ is more specific, 
		but it may include the dependence on the quadratic variation of the underlying process,
		which is related to the variance option in finance.
		Moreover, our results consider the multiple marginals case,
		such an extension of their technique seems not obvious,
		see also the work of Hou \& Ob{\l}\'oj \cite{Hou} and Biagini, Bouchard, Kardaras \& Nutz \cite{BBKN14}.
		More recently, an analogous duality is proved in the Skorokhod space under suitable conditions in Dolinsky \& Soner \cite{DS2}, where the underlying asset is assumed to take values in some subspace of c\`adl\`ag functions (see also \cite{GTT3}).

\paragraph{\bf Proof of Theorem  \ref{theo:main2}.} 
	Combining the dualities $P(\mu)=D_0(\mu)=D(\mu)$ in Theorem \ref{theo:main} and the weak dualities $\tilde P(\mu)\le \tilde D_0(\mu)\le \tilde D(\mu)$, it is enough to prove
	\b*
		P(\mu) ~~\le~~ \tilde P(\mu)
		&\mbox{and}&
		D(\mu) ~~\ge~~\tilde D(\mu),
	\e*
	where $P(\mu)$ and $D(\mu)$ are defined respectively in \eqref{eq:P} and \eqref{eq:D2} with reward function $\Phi$.

	\vspace{1mm}

	\noindent \rmi Define the process $M:=(M_t)_{0\le t\le 1}$ by 
	\b*
		M_t
		~:=~
		B_{\big(T_k+ \frac{t - t_k}{t_{k+1} -t}\big) \wedge T_{k+1}} \mbox{ for all } t \in [t_k, t_{k+1}) \mbox{ and } 0\le k \le n-1,
	\e*
with $T_0=t_0=0$ and $M_1=B_{T_n}$. It is clear that $M$ is a continuous martingale under every probability $\Pb\in\Pcb$ and $M_{t_k} = B_{T_k}$ for all $k=1,\cdots, n$, which implies in particular $M_{t_k}\stackrel{\Pb}{\sim}\mu_k$ for every $\Pb\in\Pcb(\mu)$. Let $\Pb\in\Pcb(\mu)$ be arbitrary, then $\Pt := \Pb\circ M^{-1} \in \Mc(\mu)$. Moreover, one finds  $\Pb-$a.s., $\langle M\rangle_{t_k}=T_k$ for all $k=1,\cdots, n$ and $B_t=M_{\langle M\rangle^{-1}_t}$, which yields
	\b*
		\xi(M)
		~=~
		\Phi\big( B, \langle M\rangle_{t_1}, ..., \langle M\rangle_{t_n}\big)
		~=~
		\Phi(B, T),
		~~\Pb-\mbox{a.s.}
	\e*
	Thus
	\be \label{eq:reward_equal}
		\tilde P(\mu)
		~~\ge~~
		\E^{\Pb} \big[ \xi( M) \big]
		~~=~~
		\E^{\Pb} \big[ \Phi \big(B,T\big) \big].
	\ee
	It follows that
	\b*
		P(\mu) &\le& \tilde P(\mu).
	\e*
	
	\vspace{1mm}

	\noindent \rmii Let us now prove $\tilde D(\mu) \le D(\mu)$. 
	Let $(\lambda, \Hb) \in \Dc$, i.e. $(\lambda, \Hb) \in \Lambda \x \Hcb$ be such that
	\b*
		\lambda(B_T) + (\Hb\cdot B)_{T_n}
		&\ge&
		\Phi \big(B,T \big),
		~~ \Pcb-\mbox{q.s.}.
	\e*
	For every $\tilde \P\in \Mc$, 
	it follows by Dambis-Dubins-Schwarz theorem that 
	the time-changed process $W$ defined in \eqref{eq:TimeChangeBM} is a Brownian motion with respect to the time-changed filtration $ \big( \Fct_{\langle X \rangle^{-1}_t } \big)_{t \ge 0}$
	under $\Pt$ and
	\b*
		X_t~=~W_{\Xcr_t} \mbox{ for every } t\in [0,1],~ \tilde \P-\mbox{a.s.}
	\e*
	Moreover, $\Xcr_I:=(\Xcr_{t_k})_{1\le k\le n}$ are stopping times w.r.t. the time-changed filtration $\big(\Fct_{\langle X \rangle^{-1}_t}\big)_{t\ge 0}$.
	Let us define $\Pb := \Pt\circ \big(W, \Xcr_{t_1}, ..., \Xcr_{t_n}\big)^{-1}$, 
	then $\Pb\in\Pcb$ and thus we have $\Pt-$a.s.
	\b*
		\lambda\big(W_{\Xcr_I}\big) 
		+ 
		\big(\Hb_s\cdot W\big)_{\Xcr_1}
		&\ge& \Phi \big(W, \Xcr_{t_1}, ..., \Xcr_{t_n}\big).
	\e*
	Define
	\b*
		\Ht_s(X)
		&:=&
		\Hb_{\Xcr_s} \Big(W, \Xcr_{t_1}, ..., \Xcr_{t_n}\Big),
	\e*
	then it follows by Propositions V.1.4 and V.1.5 of Revuz and Yor \cite{RevuzYor} that $\Hb$ is $\Ft-$progressively measurable such that 
	\b*
		\int_0^1 \Ht^2_s d\langle X \rangle_s
		~~=~~
		\int_0^{\langle X \rangle_1} \Hb_s^2 ds
		~~<~~ +\infty, ~\Pt -\mbox{a.s.},
	\e*
	and
	\b*
	\big(\Hb\cdot W\big)_{\Xcr_t}
		~=~
		(\tilde H\cdot X)_t \mbox{ for every } 0\le t\le 1, ~\Pt -\mbox{a.s.}
	\e*
	Hence
	\be \label{eq:superh_imme}
	\lambda(X_I) 
		+
	(\tilde H\cdot X)_1
		~~\ge~~
		\Phi\big(W, \Xcr_{t_1}, ..., \Xcr_{t_n}\big)
		~~=~~
		\xi \big(X \big),
		~ \Pt-\mbox{a.s.}	
	\ee
	Notice that $\Hb \in \Hcb$, and hence $(\Hb\cdot W)$ is a 
	strong supermartingale under $\Pt$, which implies by the time-change argument that the stochastic integral $\big(\Hb\cdot W\big)_{\Xcr_{\cdot}}$ is a supermartingale under $\Pt$ (with respect to its natural filtration) and so it is with $(\Ht\cdot X)$. Hence $\Ht \in \Hct$ and further $(\lambda, \Ht) \in \tilde \Dc$. It follows that $\tilde D(\mu) \le D(\mu)$, which concludes the proof.
	\qed

\section{Proof of Theorem \ref{theo:main} }
\label{sec:proofs}

	To prove our main results in Theorem \ref{theo:main}, 
	we start with some technical lemmas in Section \ref{subsec:technical_Lemma}.
	Then in Section \ref{sec:proof_duality1}, we provide the existence of the optimizer $\Pb^* \in \Pcb(\mu)$ and the first duality $P(\mu) = D_0(\mu)$ in Theorem \ref{theo:main},
	where the main argument is the compactness as well as the Fenchel-Moreau theorem.  
	
	Finally, in Section \ref{subsec:proof_duality2}, we complete the proofs for the second duality $P(\mu) = D(\mu)$ in Theorem \ref{theo:main} and $P(\mu) = D'(\mu) = D''(\mu)$ in Proposition \ref{prop:main_p},
	using classical results from optimal stopping theory.
	The super-hedging strategy in the second dual formulation can be obtained directly from the
	Doob-Meyer decomposition of the Snell envelop of the stopping problem in $D_0(\mu)$,
	and the martingale representation theorem.
	The argument is better illustrated in Section \ref{subsubsec:proof_prop} in the proof of Proposition \ref{prop:main_p}.

\subsection{Technical lemmas}
\label{subsec:technical_Lemma}

	Recall that $\Pbb^{\preceq}$ denotes the collection of all centered peacocks, which is a collection of vectors of probability measures on $\R$.
	We first introduce a notion of convergence $\Wc_1$ on $\Pbb^{\preceq}$ which is stronger than the weak convergence. 
	A sequence of centered peacocks $\big(\mu^m=(\mu^m_1,\cdots,\mu^m_n)\big)_{m \ge 1}\subset\Pbb^{\preceq}$ is said to converge under $\Wc_1$ to $\mu^0=(\mu^0_1,\cdots,\mu^0_n)\in\Pbb^{\preceq}$
	if $\mu_k^m$ converges to $\mu_k^0$ under the Wasserstein metric for all $k=1,\cdots, n$
	(we denote it by $\mu^m \stackrel{\Wc_1}{\lro} \mu^0$).
	Notice that the convergence under Wasserstein metric 
	is equivalent to the convergence under the weak convergence topology
	as well as the convergence of first order moment (Definition 6.1 in Villani \cite{Villani}).
	More precisely, by Theorem 6.9 of \cite{Villani}, the convergence 
	$\mu^m \stackrel{\Wc_1}{\lro} \mu^0$ holds if and only if
	\be \label{eq:Wasserstein}
		\lim_{m \to \infty}\mu^m_k(\phi) 
		~=~
		\mu^0_k(\phi) 
		~~\mbox{for all}~
		\phi\in\Cc_1
		~~\mbox{and}~~k=1,\cdots, n.
	\ee

	Further, in order to apply the Fenchel-Moreau theorem, we shall consider a linear topological space containing all centered peacocks.
	Let $\Mb$ denote the space of all finite signed measures $\nu$ on $\R$ such that
		$\int_{\R} \big(1+|x| \big) ~ |\nu|(dx)
		<
		+\infty$. 
	We endow $\Mb$ with Wasserstein topology,
	i.e. for $(\nu^m)_{m \ge 1}\subset\Mb$ and $\nu^0\in\Mb$, 
	we say $\nu^m$ converges to $\nu^0$ under $\Wc_1$ if
	\be \label{eq:Wasserstein}
		\lim_{m \to \infty} 
		\int_{\R} \phi(x) \nu^m  (dx) 
		~=~
		\int_{\R} \phi(x)  \nu^0 (dx)
		~~~\mbox{for all}~~ \phi \in \Cc_1.
	\ee
	Let $\Mb^n:=\Mb \x...\x\Mb$ be the $n-$product of $\Mb$, endowed with the product topology.
	It is clear that under $\Wc_1$, $\Pbc$ is a closed convex subspace of $\Mb^n$ and the restriction of this convergence on $\Pbc$ is the same as the Wasserstein convergence.

	It is well known that the space of all finite signed measures equipped with the weak convergence topology
	is a locally convex topological vector space,
	and its dual space is the space of all bounded continuous functions (see e.g. Section 3.2 of  Deuschel \& Stroock \cite{DS0}).
	By exactly the same arguments (see Appendix of \cite{GTT3}), we have the following similar result.

	\begin{Lemma} \label{lemm:dualM}
		There exists a topology $\Oc_n$ for $\Mb^n$ which is compatible with the $\Wc_1-$convergence,
		such that $(\Mb^n,\Oc_n)$ is a Hausdorff locally convex space. 
		Moreover, its dual space is $(\Mb^n)^{\ast}=\Lambda$.
	\end{Lemma}

	We next turn to the space $\Pcb(\Omb)$ of all Borel probability measures on the Polish space $\Omb$.
	Denote by $C_b(\Omb)$ the collection of all bounded continuous functions on $\Omb$, and $B_{mc}(\Omb)$ the collection of all bounded measurable function $\phi$, such that $\theta \mapsto \phi(\om, \theta)$ is continuous for all $\om \in \Om$.
	Notice that the weak convergence topology on $\Pcb(\Omb)$ is defined as the coarsest topology under which $\Pb \mapsto \E^{\Pb}[\xi]$ is continuous for all $\xi \in C_b(\Omb)$.
	Following Jacod \& M\'emin \cite{JacodMemin}, we introduce the stable convergence topology on $\Pcb(\Omb)$ as the coarsest topology under which $\Pb \mapsto \E^{\Pb}[\xi]$ is continuous for all $\xi \in B_{mc}(\Omb)$.
	Recall that every probability measure in $\Pcb$  (defined by \eqref{eq:def_Pcb}) has the same marginal law on $\Om$.
	Then as an immediate consequence  of Proposition 2.4 of \cite{JacodMemin}, we have the following result.
	\begin{Lemma} \label{lemm:stable_cvg}
		The weak convergence topology and the stable convergence topology 
		coincide on the space $\Pcb$.
	\end{Lemma}

	\begin{Lemma}\label{lem:tightness}
		Let $(\mu^m)_{m\ge 1}$ be a sequence of centered peacocks such that 
		$\mu^m\stackrel{\Wc_1}{\lro}\mu^0$,
		and $(\Pb_m)_{m\ge 1}$ 
		a sequence of probability measures with $\Pb_m\in\Pcb(\mu^m)$ for all $m\ge 1$. 
		Then $(\Pb_m)_{m\ge 1}$ is relatively compact under the weak convergence topology.
		Moreover, any accumulation point of $(\Pb_m)_{m\ge 1}$ belongs to $\Pcb(\mu^0)$.
	\end{Lemma}

\noindent {\bf Proof.} \rmi 
	For any $\eps>0$, there exists a compact set 
	$D\subset \Om$ such that $\Pb_m(D \x\Theta)=\P_0(D )\ge 1-\eps$ for every $m\ge 1$. In addition, by Proposition 7 of Monroe \cite{Monroe}, one has for any constant $C > 0$,
	\b*
		\Pb_m \big[T_n\ge C \big]
		~\le~
		C^{-1/3} 
		\Big(1 + \big( \mu_n^m(|x|)\big)^2\Big)~\le~C^{-1/3} 
		\Big( 1 + \big( \sup_{m\ge 1}\mu_n^m(|x|)\big)^2 \Big).
	\e*
	Choose the cube $[0,C]^n$ large enough such that $\Pb_m\big[T\in [0,C]^n\big]\ge 1-\eps$ for all $m\ge 1$. 	The tightness of $(\Pb_m)_{m \ge 1}$ under weak convergence topology follows by
	\b*
		\Pb_m\big[D\times [0,C]^n\big]~\ge ~\Pb_m\big[D\times \Theta\big]+\Pb_m\big[\Om\times [0,C]^n\big]-1~		\ge~1-2\eps \mbox{ for all } m\ge 1.
	\e*
	Let $\Pb_0$ be any limit point. By possibly subtracting a subsequence, we assume that $\Pb_m\to\Pb_0$ weakly.

	\vspace{1mm}

	\noindent \rmii Notice that $B$ is $\Fbb-$Brownian motion under each $\Pb_m$ and thus the process
$\varphi(B_t) - \int_0^t \frac{1}{2} \varphi''(B_s) ds$ is a $\Fbb-$martingale under $\Pb_m$ whenever $\varphi$ is bounded, smooth and of bounded derivatives. 
	Notice that the maps $(\om,t)\mapsto\varphi(\om_t)-\int_0^t\varphi''(\om_s)ds$ is also bounded continuous, then
	$$
		\E^{\Pb_m} 
		\Big[ 
			\Big( \varphi(B_t) - \varphi(B_r) 
			-
			\int_r^t \frac{1}{2} \varphi''(B_u) du \Big)
			~\psi
		\Big]
		~=~ 0,
 	$$
 	for every $s < r < t$ and bounded continuous and $\Fcb_r-$measurable random variable $\psi$.
	Taking the limit $m \to \infty$, it follows that 
	\be \label{eq:mart_pb_P0}
		\E^{\Pb_0} 
		\Big[ 
			\Big( \varphi(B_t) - \varphi(B_r) 
			-
			\int_r^t \frac{1}{2} \varphi''(B_u) du \Big)
			~\psi
		\Big]
		&=& 0,
 	\ee
	for all $\Fcb_r-$measurable and bounded continuous random variables $\psi$.
	Since $\Fcb_s \subset \Fcb_{r-}$, where $\Fcb_{r-}$ is generated by the class of all $\Fcb_r-$measurable bounded continuous random variables (see Lemma \ref{lemm:Fcb}),
	it follows that \eqref{eq:mart_pb_P0} is still true for every bounded and $\Fcb_s-$measurable $\psi$.
	Letting $r \to s$, by the dominated convergence theorem, 
	it follows that \eqref{eq:mart_pb_P0} holds for every $s < t$ and bounded $\Fcb_s-$measurable random variable $\psi$.
	This implies that 
	$B$ is an $\Fbb-$Brownian motion under $\Pb_0$.

	\vspace{1mm}
	
	\noindent \rmiii We next assume that $\Pb_m \in \Pcb(\mu^m)$ and prove
	\be \label{eq:claim_ui}
		B_{T_n \wedge\cdot}
		\mbox{ is uniformly integrable under } \Pb_0.
	\ee
	The convergence of $(\mu^m)_{m\ge 1}$ to $\mu^0$ implies in particular
	\b*
		\E^{\Pb_m} \big[ \big( \big| B_{T_n} \big| - R \big)^+ \big]
		~=~
		\mu_n^m\big((|x| -R)^+\big)
		~\lro~
	\mu_n^0\big((|x| -R)^+\big)
		~<~ +\infty.
	\e*
	Therefore, for every $\eps > 0$, there is $R_{\eps} > 0$ large enough such that $\mu_n^m\big((|x| -R_{\eps})^+\big)<\eps$ for every $m\ge 1$. It follows by Jensen's inequality and $|x| \1_{\{|x| > 2R\}} \le 2 (|x| -R)^+$ that
	\b*
		\E^{\Pb_m} \big[\big| B_{T_n \wedge t} \big| 
		\1_{\{| B_{T_n \wedge t} | > 2R_{\eps} \}} \big]
		~\le~
		2 \E^{\Pb_m} \big[ \big( \big| B_{T_n} \big| - R_{\eps} \big)^+ \big]
		~\le~
		2 \eps ~\mbox{ for all } t \ge 0.
	\e*
	Notice also that the function $|x|\1_{\{|x|>2R_{\eps}\}}$ is lower semicontinuous and we obtain by Fatou's lemma
	\b*
		\E^{\Pb_0} \big[\big| B_{T_n \wedge t} \big| 
		\1_{\{| B_{T_n \wedge t} | > 2R_{\eps} \}} \big]
		&\le&\liminf_{m\to\infty}\E^{\Pb_m} \big[\big| B_{T_n \wedge t} \big| 
		\1_{\{| B_{T_n \wedge t} | > 2R_{\eps} \}} \big]~~\le~~2\eps,
	\e*
	which justifies the claim \eqref{eq:claim_ui}.
	Moreover, since the map $(\om,\theta) \mapsto \om_{\theta_k}$ is continuous, 
	it follows that $B_{T_k}\stackrel{\Pb_0}{\sim} \mu^0_{k}$ for all $k = 1, \cdots, n$.
	Therefore, $\Pb_0 \in \Pcb(\mu^0)$, which concludes the proof.
	 \qed

\subsection{Proof of the first duality}
\label{sec:proof_duality1}

	We now provide the proof of the first duality result in Theorem \ref{theo:main}.
	The main idea is to show that $\mu \mapsto P(\mu)$ is concave and upper-semicontinuous
	and then to use the Fenchel-Moreau theorem.

	\begin{Lemma} \label{lemm:ucs}
		Under Assumption \ref{Hyp:phi}, the map $\mu \in \Pbb^{\preceq} \mapsto P(\mu) \in \R$ 
		is concave and upper-semicontinuous w.r.t. $\Wc_1$.
		Moreover, for every $\mu \in \Pbb^{\preceq}$, there is some $\Pb^* \in \Pcb(\mu)$ such that
		$\E^{\Pb^*}[ \Phi] = P(\mu)$.
	\end{Lemma}
	\noindent{\bf Proof.} \rmi Let $\mu^1, \mu^2 \in \Pbb^{\preceq}$, 
	$\Pb_1 \in \Pcb(\mu^1)$ and $\Pb_2 \in \Pcb(\mu^2)$ and $\alpha \in (0,1)$,
	then by their definition, one has $\alpha \Pb_1 + (1- \alpha) \Pb_2 \in \Pcb( \alpha \mu^1 + (1- \alpha) \mu^2)$.
	It follows immediately that the map $\mu \mapsto P(\mu)$ is concave.
	
	\vspace{1mm}

	\noindent \rmii We now prove that $\mu \mapsto P(\mu)$ is upper-semicontinuous w.r.t. $\Wc_1$.
	Let $(\mu^m)_{m \ge 1} \subset  \Pbb^{\preceq}$ and $\mu^m \to \mu^0 \in \Pbb^{\preceq}$ in $\Wc_1$.	
	After possibly passing to a subsequence, we can have a family  $(\Pb_m)_{m\ge1}$ such that
	\b*
		\Pb_m ~\in~ \Pcb(\mu^m)
		~~\mbox{and}~~
		\limsup_{m \to \infty} P(\mu^m)
		~=~
		\lim_{m \to \infty} \E^{\Pb_m} \Big[ \Phi \big(B,T \big) \Big].
	\e*
	By Lemma \ref{lem:tightness}, we may find a subsequence still denoted by $(\Pb_m)_{m \ge 1}$, which converges weakly to some $\Pb_0\in\Pcb(\mu^0)$.
	Moreover, it follows by Lemma \ref{lemm:stable_cvg} that
	the map $\P \mapsto \E^{\P} \big[ \Phi(B,T) \big]$ is upper-semicontinuous on $\Pcb$ w.r.t. the weak convergence topology
	for all $\Phi$ satisfying Assumption \ref{Hyp:phi}.
	We then obtain by Fatou's lemma that
	\b*
		\limsup_{m \to \infty} P(\mu^m)
		~=~
		\lim_{m \to \infty} \E^{\Pb_m} \Big[ \Phi \big(B,T \big) \Big]
		~\le~
		\E^{\Pb_0} \Big[ \Phi \big(B,T \big) \Big]
		~\le~
		P(\mu^0).
	\e*
	
	\noindent \rmiii Let $\mu \in \Pbb^{\preceq}$, choosing $\mu^m = \mu$ and using the same arguments,
	it follows immediately that there is some $\Pb^* \in \Pcb(\mu)$ such that
		$\E^{\Pb^*}[ \Phi] = P(\mu)$.
	\qed
	
	\vspace{2mm}

	\noindent The results in Lemma \ref{lemm:ucs} together with the Fenchel-Moreau theorem implies the first duality in Theorem  \ref{theo:main}.
	Before providing the proof, we consider the optimal stopping problem arising in the dual formulation \eqref{eq:D1}. Denote for every $\lambda \in \Lambda$,
	\be \label{eq:phi_lambda}
		\Phi^{\lambda}(\om,\theta)
		~:=~
		\Phi(\om,\theta)
		-\lambda(\om_{\theta}) \mbox{ for all } (\om,\theta) \in \Omb.
	\ee
	Recall that $\Tc^a$ denotes the collection of all increasing families of $\F^a-$stopping times 
	$\tau = (\tau_1, \cdots, \tau_n)$ such that $B_{\tau_n \wedge \cdot}$ is uniformly integrable.
	Recall also $\Pcb$ is defined in \eqref{eq:def_Pcb} as set of measures of the Brownian motion and stopping times.
	Let $N > 0$, denote also by $\Tc^a_N \subset \Tc^a$ the subset of families 
	$\tau = (\tau_1, \cdots, \tau_n)$ such that $\tau_n \le N$, $\P_0-$a.s.
	Denote further by $\Pcb_N \subset \Pcb$ the collection of $\Pb \in \Pcb$ such that $T_n \le N$, $\Pb-$a.s.

	\begin{Lemma} \label{lemm:equiv_OS}
		Let $\Phi$ be bounded, then for every $\lambda \in \Lambda$,
		\be \label{eq:equiv_OS}
			\sup_{\tau \in \Tc^a} \E^{\P_0}[ \Phi^{\lambda}(B,\tau) ]
			&=& 
			\lim_{N \to \infty} \sup_{\tau \in \Tc^a_N} \E^{\P_0}[ \Phi^{\lambda}(B,\tau) ] \\
			&=&
			\lim_{N \to \infty} \sup_{\Pb \in \Pcb_N} 
			\E^{\Pb} \big[ \Phi^{\lambda}(B,T) \big]
			~=~
			\sup_{\Pb \in \Pcb} \E^{\Pb} \big[ \Phi^{\lambda}(B,T) \big]. \nonumber
		\ee
		In particular, let $\phi \in \Cc_1$ and denote by $\phi^{conc}$ its concave envelope, one has
		\be \label{eq:phi_conc}
			\sup_{\tau \in \Tc^a} \E^{\P_0} [\phi (B_{\tau_n}) ]
			&=&
			\phi^{conc}(0).
		\ee
	\end{Lemma}
	\noindent{\bf Proof.} \rmi Given $\lambda \in \Lambda$, there is some constant $C > 0$ such that 
	\be \label{eq:Phi_lambda}
		\big| \Phi^{\lambda}(B,\tau) \big| 
		~~\le~~
		C \Big(1 + \sum_{k =1}^n \big| B_{\tau_k} \big| \Big).
	\ee
	Let $\tau \in \Tc^a$, define $\tau^N := (\tau^N_1, \cdots, \tau^N_n)$ with $\tau^N_k := \tau_k \wedge N$,
	then it is clear that $\lim_{N \to \infty} \Phi^{\lambda}(B,\tau^N) = \Phi^{\lambda}(B,\tau)$, $\P_0-$a.s.
	By the domination in \eqref{eq:Phi_lambda} and the fact that $B_{\tau_n \wedge \cdot}$ is uniformly integrable,
	we have $\lim_{N \to \infty} \E^{\P_0}\big[ \Phi^{\lambda}( B,\tau^N) \big]= \E^{\P_0} \big[\Phi^{\lambda}(B,\tau) \big]$.
	It follows by the arbitrariness of $\tau \in \Tc^a$ and the fact $\Tc^a_N \subset \Tc^a$ that
	\b*
		\sup_{\tau \in \Tc^a} \E^{\P_0}[ \Phi^{\lambda}(B,\tau) ]
		&=& 
		\lim_{N \to \infty} \sup_{\tau \in \Tc^a_N} \E^{\P_0}[ \Phi^{\lambda}(B,\tau) ].		
	\e*
	By the same arguments, it is clear that we also have
	\b*
		\sup_{\Pb \in \Pcb} \E^{\Pb} \big[ \Phi^{\lambda}(B,T) \big]
		&=&
		\lim_{N \to \infty} \sup_{\Pb \in \Pcb_N} 
		\E^{\Pb} \big[ \Phi^{\lambda}(B,T) \big].
	\e*

	\noindent \rmii We now apply Lemma \ref{lemm:OS_equiv}
	to prove that for every fixed constant $N > 0$,
	\be \label{eq:equiv_N}
		\sup_{\tau \in \Tc^a_N} \E^{\P_0}[ \Phi^{\lambda}(B,\tau) ]
		&=&
		\sup_{\Pb \in \Pcb_N} \E^{\Pb} \big[ \Phi^{\lambda}(B,T) \big].
	\ee
	First, let us suppose that $n = 1$.
	Let $\Pb \in \Pcb_N$, denote $Y_t := \Phi^{\lambda}(B,t \wedge N)$,
	it is clear that $\E^{\Pb} \big[ \sup_{t \ge 0} Y_t \big] < \infty$.
	Denote by $\Fbb^{\Pb} = (\Fcb^{\Pb}_t)_{t \ge 0}$ the augmented filtration of $\Fbb$ under $\Pb$ and by
	$\Fbb^{B,\Pb}$ the filtration generated by $B$ on $\Omb$ and by
	$\Fbb^{B,\Pb} = (\Fcb^{B,\Pb}_t)_{t \ge 0}$ its $\Pb-$augmented filtration.
	It is clear that $\Fcb^{B, \Pb}_t \subset \Fcb_t^{\Pb}$.
	More importantly, by the fact that $B$ is a $\Fbb^{\Pb}-$Brownian motion under $\Pb$,
	it is easy to check that the probability space $(\Omb, \Fcb^{\Pb}, \Pb)$ together with the filtration $\Fbb^{\Pb}$ and $\Fbb^{B,\Pb}$ satisfies Hypothesis (K) (Assumption \ref{hyp:K}).
	Then by Lemma \ref{lemm:OS_equiv},
	$\E^{\Pb} \big[ \Phi^{\lambda}(B,T) \big] 
	\le
	\sup_{\tau \in \Tc^a_N} \E^{\P_0}[ \Phi^{\lambda}(B,\tau) ]$
	and hence
	$\sup_{\Pb \in \Pcb_N} \E^{\Pb} \big[ \Phi^{\lambda}(B,T) \big] 
	\le
	\sup_{\tau \in \Tc^a_N} \E^{\P_0}[ \Phi^{\lambda}(B,\tau) ]$.
	We then have equality \eqref{eq:equiv_N}
	since the inverse inequality is clear.
	Finally, when $n > 1$, it is enough to use the same arguments together with induction to prove \eqref{eq:equiv_N}.

	\vspace{1mm}

	\noindent \rmiii To prove \eqref{eq:phi_conc} it suffices to set $\Phi \equiv 0$ and $n = 1$.
	Then by \eqref{eq:equiv_OS}, it follows that
	\b*
		\sup_{\tau \in \Tc^a} \E^{\P_0} \big[ \phi(B_{\tau}) \big]
		~~=~~
		\lim_{N \to \infty}
		\sup_{\tau \in \Tc^a_N} \E^{\P_0} \big[ \phi(B_{\tau}) \big]
		~~\le~~ \phi^{conc}(0).
	\e*
	The inverse inequality is obvious by considering the exiting time of the Brownian motion from an open interval.
	We hence conclude the proof. \qed

	\vspace{2mm}

	\noindent {\bf Proof of Theorem \ref{theo:main} $\mathrm{(i)}$.}
	The existence of optimal embedding is already proved in Lemma \ref{lemm:ucs}.
	For the first duality result, we shall use the Fenchel-Moreau theorem.
	Let us first extend the map $ \mu \mapsto P(\mu)$ from $\Pbb^{\preceq}$ to $\M^n$ 
	by setting that $P(\mu) = -\infty$,
	for every $\mu \in \M^n \setminus \Pbb^{\preceq}$.
	It is easy to check, using Lemma \ref{lemm:ucs}, 
	that the extended map $\mu \mapsto P(\mu)$ from the topological vector space $\M^n$ to $\R$ is still concave and upper-semicontinuous.
	Then by Fenchel-Moreau theorem together with Lemma \ref{lemm:dualM}, it follows that
	\b*
		P(\mu)~=~P^{\ast\ast}(\mu)
		&=&\inf_{\lambda\in\Lambda} 
		\Big\{
			\sup_{\nu\in\Pbb^{\preceq}}
			\sup_{\Pb\in\Pcb(\nu)}
			\E^{\Pb} 
			\big[ \Phi^{\lambda}(B,T) \big]
			+\mu(\lambda)
		\Big\} \\
		&=& 
		\inf_{\lambda \in \Lambda} \Big\{
			\sup_{\tau \in \Tc^a} \E^{\P_0}
			\big[
				\Phi ^{\lambda}\big(B,\tau\big) 
			\big] 
			+ \mu(\lambda)
		\Big\},
\e*
	where the last equality follows by \eqref{eq:equiv_OS}. 
	Hence we have $P(\mu) = D_0(\mu)$.
	\qed

	\begin{Remark} \label{rem:lambda_positive}
		When $\Phi$ is bounded (which is the relevant case by the reduction of Section \ref{subsubsec:reduction}), we can prove further that
		\be \label{eq:D1_lambda_plus}
			D_0(\mu)
			&=&
			\inf_{\lambda \in \Lambda^+} \Big\{
				\sup_{\tau \in \Tc^a} \E^{\P_0}
				\big[
					\Phi^{\lambda} (B,\tau) 
				\big] 
				+ \mu(\lambda)
			\Big\},
		\ee
		where
		\b*
			\Lambda^+
			&:=&
			\big\{ \lambda = (\lambda_1, \cdots, \lambda_n) \in \Lambda ~: 
			\lambda_k \ge 0 \mbox{ for all } k =1, \cdots, n
			\big\}.
		\e*
		Indeed, using \eqref{eq:phi_conc}, it is easy to see that in the definition of $D_0(\mu)$, it is enough to take the infimum over the class of all functions $\lambda \in \Lambda^+$ such that the convex envelope $\lambda_k^{conv}(0)>-\infty$ for all $k=1, \cdots, m$,
		since by \eqref{eq:phi_conc} and the boundedness of $\Phi$,
		$\sup_{\tau\in\Tc^a}\E^{\P_0}[\Phi^{\lambda}]=+\infty$ whenver $(-\lambda_k)^{conc}(0)= \infty$ for some $k$. 
		Hence the infimum is taken among all $\lambda\in\Lambda$ such that $\lambda_k^{conv}(0) > -\infty$ for all $k=1,\cdots, m$,
		and consequently $\lambda_k$ is dominated from below by some affine function. 
		Since $\E^{\P_0} [B_{\tau_k}] = 0$ for every $\tau \in \Tc^a$, 
		we see that by possibly subtracting from $\lambda_k$ the last affine function, 
		it is enough to take infimum over the class $\Lambda^+$.
	\end{Remark}

\subsection{Proof of the second duality}
\label{subsec:proof_duality2}

	We now prove the second duality $P(\mu) = D(\mu)$ in Theorem \ref{theo:main} $\mathrm{(ii)}$, 
	and $P(\mu) = D'(\mu) = D''(\mu)$ in Proposition \ref{prop:main_p}.
	The main technique is to use the Snell envelope characterization of the optimal stopping problem, together with the Doob-Meyer decomposition.
	We will provide the proof progressively.
	In Section \ref{subsubsec:weakduality}, we prove a weak duality result $P(\mu) \le D(\mu)$ in the context of Theorem \ref{theo:main} and Proposition \ref{prop:main_p}.
	Then in Section \ref{subsubsec:reduction},
	we show that it is enough to prove Theorem \ref{theo:main} \rmii and Proposition \ref{prop:main_p}
	for bounded reward function $\Phi$.
	Next, in Section \ref{subsubsec:proof_prop}, we provide the proof of Proposition \ref{prop:main_p},
	which implies immediately Theorem \ref{theo:main} \rmii under Assumptions \ref{Hyp:phi} and \ref{Hyp:phi2} $\mathrm{(i)}$.
	Finally, we complete the proof of Theorem \ref{theo:main} \rmii 
	under Assumptions  \ref{Hyp:phi} and \ref{Hyp:phi2} \rmii or \rmiii
	in Sections \ref{subsubsec:Duality2_1} and \ref{subsubsec:duality2usc}.

	Throughout this subsection,
	we say a process $X$, on filtered space $(\Om, \Fc, \P_0, \F^a)$, is of class (DL) if for each $t \ge 0$, 
	the family $\{ X_{\tau} ~:\tau \in \Tc^a, \tau \le t\}$ is uniformly integrable;
	we say an $\F^a-$optional process $X$ of class (DL) is a supermartingale if for all bounded stopping times $\sigma \le \tau$, one has $X_{\sigma} \ge \E^{\P_0}[ X_{\tau} | \Fc^a_{\sigma}]$.

\subsubsection{On the weak duality}
\label{subsubsec:weakduality}

	We notice that from their definition, we can easily have the weak duality
	in the context of Theorem \ref{theo:main} and Proposition \ref{prop:main_p}.

	\begin{Lemma} \label{lemm:weakduality}
		Let $\Phi : \Omb \to \R$ be non-anticipative, then one has
		\b*
			P(\mu) ~\le~ D(\mu),
			~~~~
			P(\mu) ~\le~ D'(\mu)
			~~~\mbox{and}~~
			P(\mu) ~\le~ D''(\mu).
		\e*
	\end{Lemma}
	\proof 
		\rmi Let $(\lambda, \Hb) \in \Dc$,
		one has, by its definition,
		\b*
			\lambda(B_T) + (\Hb \cdot B)_{T_n}~ \ge~ \Phi(B,T), ~~\Pcb- \mbox{q.s.},
		\e*
		where $(\Hb \cdot B)$ is a strong supermartingale.
		Let $\Pb \in \Pcb(\mu)$, taking the expectation of the above inequality under $\Pb$,
		it follows that 
		\b*
			\mu(\lambda) ~=~ \E^{\Pb}[ \lambda(B_T) ]  ~\ge~ \E^{\Pb} \big[ \Phi(B,T) \big].
		\e*
		It follows by the arbitrariness of $\Pb \in \Pcb(\mu)$ and $(\lambda, \Hb) \in \Dc$ that one has $P(\mu) \le D(\mu)$.
		
		\vspace{1mm}
		
		\noindent \rmii We next prove $P(\mu) \le D'(\mu)$.
		Notice that for any $H^0 \in \Hc$, one has $(H^0 \cdot B)_t \ge - C(1 + |B_t|)$.
		Then for any $\Pb \in \Pcb$, since the process $B_{T_n \wedge \cdot}$ is uniformly integrable,
		it follows by Fatou's lemma that 
		\b*
			\E^{\Pb} \Big[  \int_{T_k}^{T_{k+1}} 
				H^0_s d B_s
			\Big]
			~\le~ 0,
			~~~\mbox{for each}~k=0, \cdots, n-1.
		\e*
		Using exactly the same arguments as above, we can conclude that $P(\mu) \le D'(\mu)$.
		
		\vspace{1mm}
		
		\noindent \rmiii Similarly, one can easily prove that $P(\mu) \le D''(\mu)$.
	\qed

\subsubsection{Reduction to bounded reward functions}
\label{subsubsec:reduction}

	\begin{Proposition} \label{prop:Phi_bounded}
		To prove Theorem \ref{theo:main} $\mathrm{(ii)}$ and Proposition \ref{prop:main_p},
		it is enough to prove the results
		under additional condition that $\Phi$ is bounded.
	\end{Proposition}
	\noindent{\bf Proof.} We will prove it in the context of  Theorem \ref{theo:main} $\mathrm{(ii)}$,
	since the arguments in the context of Proposition \ref{prop:main_p} is the same.

	Assume that the duality $P(\mu) = D(\mu)$ holds true
	whenever $\Phi$ is bounded and satisfies Assumptions \ref{Hyp:phi} and \ref{Hyp:phi2}.

	We now consider the case without boundedness of $\Phi$.
	Let $\Phi_m := \Phi \vee (-m)$ (or $\Phi_m := \sum_{k=1}^n (-m) \vee  \Phi_k $ in case of Assumption \ref{Hyp:phi2} $\mathrm{(iii)}$),
	then $\Phi_m$ is bounded and satisfies Assumptions \ref{Hyp:phi} and \ref{Hyp:phi2}.
	Denote by $P^m(\mu)$ and $D^m(\mu)$ the corresponding primal and dual values associated to the reward function $\Phi_m$, so that we have the duality
	\b*
		P^m(\mu)
		&=&
		D^m(\mu).
	\e*
	Further, notice that $\Phi_m\ge\Phi$, one has $P^m(\mu)=D^m(\mu)\ge D(\mu)\ge P(\mu)$,
	where the last inequality is the weak duality in Lemma \ref{lemm:weakduality}.
	Then it is enough to show that
	\b*
		\limsup_{m\to\infty}P^m(\mu)
		&\le&
		P(\mu).
	\e*
	Let $\Pb_m\in\Pcb(\mu)$ such that $\limsup_{m\to\infty}P^m(\mu)=\limsup_{m\to\infty}\E^{\Pb_m}[\Phi_m]$. Then after possibly passing to a subsequence we may assume that $\limsup_{m\to\infty}P^m(\mu)=\lim_{m\to\infty}\E^{\Pb_m}[\Phi_m]$.
	By Lemma \ref{lem:tightness}, we know that $(\Pb_m)_{m\ge 1}$ is tight and every limit point belongs to $\Pcb(\mu)$. Let $\Pb_0$ be a limit point of $(\Pb_m)_{m\ge 1}$, and label again the convergent subsequence by $m$, i.e. $\Pb_m\to\Pb_0$. Then by the monotone convergence theorem
	\b*
		P(\mu) ~~\ge~~ \E^{\Pb_0}[\Phi]
		&=&
		\lim_{m\to\infty}\E^{\Pb_0}[\Phi_m]
		~~=~~
		\lim_{m\to\infty}\Big(\lim_{l\to\infty}\E^{\Pb_l}[\Phi_m]\Big)\\
		&\ge&
		\lim_{m\to\infty}\Big(\lim_{l\to\infty}\E^{\Pb_l}[\Phi_l]\Big)
		~~=~~
		\limsup_{l\to\infty}P^l(\mu),
	\e*
	which is the required result.
	\qed

\subsubsection{Proof of Proposition \ref{prop:main_p}}
\label{subsubsec:proof_prop}

	By Proposition \ref{prop:Phi_bounded}, we can assume in addition that $\Phi$ is bounded without loss of generality.
	Then given the first duality $P(\mu) = D_0(\mu)$, it suffices to study the optimal stopping problem
	\be \label{eq:OStopping_n1}
		\sup_{\tau \in \Tc^a} \Big[ \Phi^{\lambda} \big(B, \tau\big) \Big]
		&=&
		\lim_{N \to \infty} \sup_{\tau \in \Tc^a_N} \E \big[ \Phi^{\lambda}(B, \tau) \big],
	\ee
	for a given $\lambda \in \Lambda^+$ (Remark \ref{rem:lambda_positive}) 
	and bounded $\Phi_k$.
	Notice that in this case, there is some $C$ such that
	\be \label{eq:DominationPhiLambda}
		- C \Big( 1 + \sum_{k=1}^n |\om_{\theta_k}| \Big) 
		~\le~
		\Phi^{\lambda}(\omb) 
		~\le~
		C.
	\ee

	Suppose that $n = 1$,
	then by Lemma \ref{lemm:Snell}, there is an $\F^a-$optional l\`adl\`ag process $(Z^{1,N}_t)$,
	for every $N \in \N$, which is an $\F^a-$supermartingale and the Snell envelope of the optimal stopping problem 
	$\sup_{\tau \in \Tc^a_N} \E[ \Phi^{\lambda}(B,\tau) ]$.
	Clearly, $Z^{1,N}$ increases in $N$.
	Moreover, since $Z^{1,N}$ dominates $\Phi^{\lambda}$, 
	by \eqref{eq:DominationPhiLambda}, 
	one has $- C(1 + |B_t|) \le Z^{1,N}_t \le C$ for some constant $C$ independent of $N$.
	Then by the dominated convergence theorem together with Lemma \ref{lemm:equiv_OS}, $Z^1 := \sup_{N \in \N} Z^{1,N}$ is still a l\`adl\`ag $\F^a-$supermartingale, of class (DL), such that
	$$
		Z^1_0 = \sup_{\tau \in \Tc^a} \E[ \Phi^{\lambda}(B, \tau) ],
		~\mbox{and}~~
		Z^1_t \ge \Phi^{\lambda}(B, t), ~\mbox{for all}~ t \ge 0, ~\P_0-\mbox{a.s.}
	$$
		
	Now, by the Doob-Meyer decomposition (see e.g. Lemma \ref{lemm:Mertens} below)
	for supermartingales of class (DL) without right-continuity,
	together with the martingale representation theorem,
	there is an $\F^a-$predictable process $H^1$ such that
	\b*
		\lambda_1 (B_t) ~+~ (H^1 \cdot B)_t
		~\ge~
		\Phi(B, t),
		~\mbox{for all}~ t \ge 0,
		&\P_0 -\mbox{a.s.}
	\e*
	Further, since any $\F^a-$predictable process (or equivalently $\F^a-$optional process) is indistinguishable to an $\F-$predictable process
	(see e.g. Theorem IV.78 and Remark IV.74 of Dellacherie \& Meyer \cite{DM}),
	we can also choose $H^1$ to be $\F-$predictable.
	This proves in particular that
	\b*
		D'(\mu)
		~\le~
		D_0(\mu)
		~=~
		P(\mu)
		~~\mbox{and}~~
		D''(\mu)
		~\le~
		D_0(\mu)
		~=~
		P(\mu).
	\e*
	Combining with the weak duality $P(\mu)\le D'(\mu)$ and $P(\mu) \le D''(\mu)$
	in Lemma \ref{lemm:weakduality}, we obtain
	\b*
		P(\mu)
		&=&
		D_0(\mu)
		~=~
		D'(\mu)
		~=~
		D''(\mu).
	\e*
	
	Suppose now $n = 2$, we first consider the optimal stopping problem
	\b*
		\sup_{\tau \in \Tc^0_N} \big[ \Phi_2 (B, \tau) - \lambda_2 (B_{\tau}) \big],
	\e*
	whose Snell envelope is given by $Z^{2,N}$ by Lemma \ref{lemm:Snell}, 
	where in particular $ - C(1 + |B_t|) \le Z^{2,N}_t \le C$ for some constant $C$ independent of $N$,
	and
	\b*
		Z^{2,N}_{\theta_2} ~\ge~ \Phi_2(B, \theta_2) - \lambda_2(B_{\theta_2}),~~\mbox{for all}~\theta_2 \le N, ~~\P_0- \mbox{a.s}.
	\e*
	We then reduce the multiple optimal stopping problem \eqref{eq:OStopping_n1} to the $n=1$ case, i.e.
	\b*
		\sup_{\tau \in \Tc^a_N} \E \big[ \Phi^{\lambda}(B, \tau) \big]
		~=~
		\sup_{\tau_1 \in \Tc^0_N} \E \big[ Z^{2,N}_{\tau_1}  + \Phi_1(B, \tau_1) - \lambda_1(B_{\tau_1}) \big].
	\e*
	Using again the procedure for the case $n=1$, we obtain a new Snell envelop, denoted by $Z^{1,N}$,
	such that $Z^{1,N}_t \ge - C(1 + |B_t|)$.
	
	Thus, $Z^{1,N}, Z^{2,N}$ are both supermartingales of class (D), bounded from above by $C$, and dominated from below by $-C(1 + |B_t|)$ for some constant $C > 0$ independent of $N$.
	More importantly, we have $Z^{1,N}_0 =\sup_{\tau \in \Tc^a_N} \E \big[ \Phi^{\lambda}(B, \tau) \big]$, and
	\b*
		Z^{1,N}_{\theta_1} ~+~ \big( Z^{2,N}_{\theta_2} - Z^{2,N}_{\theta_1} \big)
		&\ge&
		\Phi^{\lambda}(B, \theta_1, \theta_2),
		~\mbox{for all}~\theta_1 \le \theta_2 \le N,~~\P_0-\mbox{a.s.}
	\e*
	Since $Z^{1,N}$ and $Z^{2,N}$ both increase in $N$, define $Z^1 := \sup_{N} Z^{1,N}$ and $Z^2 := \sup_N Z^{2,N}$, it follows by the dominated convergence theorem that $Z^1$ and $Z^2$ are both supermartingales of class (DL).
	Moreover, it follows from Lemma \ref{lemm:equiv_OS} that
	$Z^1_0 = \sup_{\tau \in \Tc^a} \E \big[ \Phi^{\lambda}(B, \tau) \big]$ and
	\b*
		Z^1_{\theta_1} ~+~ \big( Z^2_{\theta_2} - Z^2_{\theta_1} \big)
		&\ge&
		\Phi^{\lambda}(B, \theta_1, \theta_2),
		~\mbox{for all}~\theta_1 \le \theta_2 ,~~\P_0-\mbox{a.s.}
	\e*
	Then $(S^1, S^2) := (Z^1, Z^2)$ are the required supermartingale in dual formulation $\Dc''$.
	Further, using the Doob-Meyer decomposition, together with the martingale representation on $Z^1$ and $Z^2$, 
	we obtain the process
	$H = (H^1, H^2)$ as we need in the dual formulation $\Dc'$.
	
	Finally, the case $n > 2$ can be handled by exactly the same recursive arguments as for the case $n =2$.
	\qed

\subsubsection{Proof of Theorem \ref{theo:main} \rmii under Assumption \ref{Hyp:phi2} \rmi }

	When $n=1$, Theorem \ref{theo:main} is an immediate consequence of Proposition \ref{prop:main_p}.

\subsubsection{Proof of Theorem \ref{theo:main} \rmii under Assumption \ref{Hyp:phi2} \rmiii }
\label{subsubsec:Duality2_1}

	Let $N >0$, we first study the multiple optimal stopping problem
	\be \label{eq:OStopping_n2}
		\sup_{\tau \in \Tc^a_N}
		\E^{\P_0} \big[ \Phi^{\lambda}(B, \tau) \big]
		~=~
		\sup_{\tau \in \Tc^a_N}
		~ \E^{\P_0} \Big[ \sum_{k=1}^n \Big( 
			\Phi_k (B, \tau_1,\cdots, \tau_k ) - \lambda_k (B_{\tau_k}) 
		\Big) 
		\Big ],
	\ee
	where $\lambda \in \Lambda^+$ and $\Phi_k$ is bounded, so that 
	\be \label{eq:Boundedness_Phi_lambda}
		- C \Big( 1 + \sum_{k=1}^n |\om_{\theta_k}| \Big) 
		~\le~
		\Phi^{\lambda}(\omb) 
		~\le~
		C,
	\ee
	for some constant $C$.
	Denote $v^N_{n+1}(\om, \theta_1, \cdots, \theta_n, \theta_n) := \Phi^{\lambda}(\om, \theta_1, \cdots, \theta_n) $.
	
	\begin{Lemma} \label{lemm:supermart_agg}
		There are functionals $(v_k^N)_{k = 1, \cdots, n}$, where $v_k^N: \Om \x (\R_+)^k \to \R$,
		such that
		$$
			v_1^N(\om, 0) ~=~ 
			\sup_{\tau \in \Tc^a_N}
			\E^{\P_0} \big[ \Phi^{\lambda}(B, \tau) \big],
		$$
		and under $\P_0$, for each $k =1, \cdots, n$, and $\theta_1 \le \cdots \le \theta_{k-1}$,
		the process 
		$$
			\theta ~\mapsto~ v_k^N(B, \theta_1, \cdots, \theta_{k-1}, \theta)
			~\mbox{is an}~ \F^a-\mbox{supermartingale},
		$$
		$$
			v^N_k(B, \theta_1, \theta_{k-1}, \theta) 
			~\ge~
			v_{k+1}^N(B, \theta_1, \cdots, \theta_{k-1}, \theta, \theta), ~\P_0-\mbox{a.s.}
		$$
		Moreover, $v^N_k$ increases in $N$ and satisfies 
		$- C(1 + \sum_{i=1}^k |\om_{\theta_i}| ) \le v^N_k(\om, \theta_1, \cdots, \theta_k) \le C $ for some constant $C$ independent of $N$.
	\end{Lemma}

	\noindent {\bf Proof of Theorem \ref{theo:main} \rmii}
	By Remark  \ref{rem:lambda_positive} and Proposition \ref{prop:Phi_bounded},
	we can assume without loss of generality that
	each $\Phi_k$ is bounded and choose $\lambda \in \Lambda^+$ in the dual formulation $D_0(\mu)$.
	
	\vspace{1mm}
	
	\noindent \rmi Let $v_k^N$ be given by Lemma \ref{lemm:supermart_agg}, we define further
	\b*
		v_k(\cdot) ~:=~ \sup_{N} v_k^N(\cdot),
		~~\mbox{so that}~
		v_1(\om,0) ~=~ \sup_{\tau \in \Tc^a}
		\E^{\P_0} \big[ \Phi^{\lambda}(B, \tau) \big].
	\e*
	It follows from the dominated convergence theorem that,
	for all  $k=1,\cdots, n$  and $0 =: \theta_0 \le \theta_1\le \cdots\le \theta_{k-1}$,
	the process $\big(v_k(B, \theta_1, \cdots, \theta_{k-1}, t )\big)_{t \ge \theta_{k-1} }$ is an $\F^a-$supermartingale
	and
	$$
			v_k(B, \theta_1, \cdots, \theta_{k-1}, t) 
			~\ge~
			v_{k+1}(B, \theta_1, \cdots, \theta_{k-1}, t, t), \mbox{for all}~t \ge \theta_{k-1},~\P_0-\mbox{a.s.}	
	$$

	\vspace{1mm}
	
	\noindent \rmii By the Doob-Meyer decomposition (see Lemma \ref{lemm:Mertens} below) and the martingale representation theorem, 
	it follows that for each $k= 1, \cdots, n$,
	there is some $\F^a-$predictable process $H^k_t(\om):=H^k_t(\om, \theta_1, ..., \theta_{k-1})$ such that
	\be \label{eq:determ_supermart}
		v_k(\om, \theta_1,..., \theta_{k-1}, \theta_{k-1}) 
		+\!\!
		\int_{\theta_{k-1}}^{\theta_{k}} H_u^k dB_u
		\!\!\!&\ge&\!\!\!
		v_k(\om, \theta_1,... \theta_{k-1}, \theta_{k}) \nonumber \\
		\!\!\!&\ge&\!\!\!
		v_{k+1}(\om, \theta_1,... \theta_{k-1}, \theta_{k}, \theta_{k}),~\P_0-\mbox{a.s.} ~~~~~~~~
	\ee
	
	\noindent \rmiii
	Next, following the pathwise construction in \eqref{eq:const_quad_cova} of the quadratic co-variation $Q^-$ of a supermartingale and a continuous martingale,
	one has a Borel version of the quadratic co-variation $\langle v_k(B, \theta_1, \cdots, \theta_{k-1}, \cdot) , B_{\cdot} \rangle_t$.
	Then by Lemma \ref{eq:pred_co_variat}, the process $H^k$ defined below is $\F-$predictable,
	\b*
		&&
		H^k_t(\theta_1, \cdots, \theta_{k-1}) \\
		&:=&
		\limsup_{\eps \to 0}
		\frac{\langle v_k(B, \theta_1, \cdots, \theta_{k-1}, \cdot) , B \rangle_t - \langle v_k (B, \theta_1, \cdots, \theta_{k-1}, \cdot) , B \rangle_{t-\eps}}{\eps}.
	\e*
	In particular, the map $(\om, \theta_1, \cdots, \theta_k) \mapsto H^k_{\theta_k} (\om, \theta_1, \cdots, \theta_{k-1})$ is Borel measurable,
	and 
	\be \label{eq:L2_loc}
		\int_{\theta_{k-1}}^{t} \big( H^k_s(\cdot, \theta_1, \cdots, \theta_{k-1}) \big)^2 ds 
		~<~
		+\infty \mbox{ for all}~ t \ge \theta_{k-1}, 
		~\P_0-\mbox{a.s.} 
	\ee

	\noindent \rmiv Next, we define a process $\Hb: \R_+\x\Omb  \to \R$ by
	\b*
		\Hb_u(\omb) 
		&:=&
		\sum_{k=1}^n \1_{(\theta_{k-1},\theta_k]} (u)
		H_u^k(\om, \theta_1, ..., \theta_{k-1}) \mbox{ for all } \omb = (\om,\theta)\in\Omb,
	\e*
	where by convention $\theta_0=0$.
	Moreover, since
	$$
		(\om, \theta_1, \cdots, \theta_k) 
		~\mapsto~
		H^k_{\theta_k} (\om, \theta_1, \cdots, \theta_{k-1}) 
		~\mbox{is Borel measurable,}~
	$$
	and one has clearly that $H^k_{\theta_k}(\om, \theta_1, \cdots, \theta_{k-1}) = H^k_{\theta_k}(\om_{\theta_k \wedge \cdot}, \theta_1, \cdots, \theta_{k-1})$,
	then the process $\Hb$ is $\Fbb-$optional by Lemma \ref{lemm:rcpd} in Appendix.

	\vspace{1mm}

	\noindent \rmv Now, let us take an arbitrary $\Pb\in\Pcb$ 
	and consider a family of r.c.p.d. (regular conditional probability distributions) $(\Pb_{\omb})_{\omb \in \Omb}$ 
	of $\Pb$ with respect to $\Fcb_{T_{k}}$ for $0\le k\le n-1$ 
	(see Lemma \ref{lemm:rcpd} for the existence of r.c.p.d.). 
	Then for  $\Pb-$almost every $\omb \in \Omb$, under the conditional probability $\Pb_{\omb}$,
	the process $t \mapsto B_t$ for $t\ge T_k$ is still a Brownian motion.
	Moreover, we have $\Pb_{\omb}( T_k = \theta_k, B_{T_k \wedge \cdot} = \om_{\theta_k \wedge \cdot}) = 1$.
	Then it follows by \eqref{eq:determ_supermart} that
	\b*
		v_{k+1}(B, T_1, ..., T_k, T_k)
		&\le&
		v_{k}(B, T_1, ..., T_k)\\
		&\le&
		v_{k}(B, T_1, ..., T_{k-1}, T_{k-1}) + \int_{T_{k-1}}^{T_{k}}H^k_s dB_s,
		~\Pb_{\omb}-\mbox{a.s.}
	\e*
	This means that the set $A_k:=\big\{v_{k+1}\le v_k+\int_{T_{k-1}}^{T_{k}}H^k_sdB_s\big\}$ is of full measure under $\Pb_{\omb}$ for $\Pb-$almost every $\omb\in\Omb$, 
	and hence by the tower property $\Pb(A_k)=1$ for all $k=0,\cdots, n$
	(we also refer to \cite{CTT} for some some discussion on the measurability of $A_k$ under $\Pb_{\omb}$).
	This yields that
	\be \label{eq:Phi_le_lambda_h}
		\Phi^{\lambda}(B, T)
		~=~ 
		v_{n+1}(B, T_1, ..., T_n, T_n)
		~\le~
		v_1(B,0)
		+
		(\Hb\cdot B)_{T_n},
		~\Pb-\mbox{a.s.}
	\ee
	
	\noindent \rmvi To conclude the proof, it suffices to check that $\Hb \in \Hcb$.
	First, for any probability measure $\Pb \in \Pcb$, by taking r.c.p.d and using \eqref{eq:L2_loc},
	it is clear that
	\b*
		\int_0^t \Hb_s^2 ds ~< +\infty \mbox{ for every } t \ge 0,~ \Pb-\mbox{a.s.}
	\e*
	Notice also that \eqref{eq:Phi_le_lambda_h} holds true for every $\Pb \in \Pcb$, and 
	by the tower property, it follows that for any $\Fbb-$stopping time $\tau$, we have for all $\Pb \in \Pcb$,
	\b*
		(\Hb\cdot B)_{T_n\wedge \tau} &\ge& - C \Big(1+ \sup_{1\le k\le n} |B_{T_k\wedge \tau}| \Big), 
		~ \Pb -\mbox{a.s.},
	\e*
	where the r.h.s. is uniformly integrable under $\Pb$.
	Using Fatou's Lemma, it follows that $(\Hb\cdot B)_{T_n\wedge \cdot}$ 
	is a strong supermartingale under every $\Pb \in \Pcb$.
	\qed

	\vspace{2mm}
	
	\noindent {\bf Proof of Lemma \ref{lemm:supermart_agg}.}
	We provide here a proof for the case $n=2$ for ease of presentation. 
	The general case can be treated by exactly the same backward iterative procedure.
	We will use the aggregation procedure in the optimal stopping theory 
	(see e.g. El Karoui \cite{EK}, Peskir \& Shiryaev \cite{PS}, Karatzas \& Shreve \cite{KS} or Kobylanski, Quenez \& Rouy-Mironescu \cite{Quenez}, etc.)
	
	{\bf 1.} For every $\tau_1 \in \Tc^0_N$, we first consider the optimal stopping problem
	\b*
		\sup_{\tau_2 \in \Tc^0_N, ~\tau_2 \ge \tau_1} 
		~\E^{\P_0} \Big[ \Phi_2 \big(B, \tau_1, \tau_2) - \lambda_2(B_{\tau_2}) \Big],
	\e*
	whose Snell envelope is denoted by $(Z^{2,N}_{\tau_1, t})_{\tau_1 \le t \le N}$.
	We shall prove in Step 2 below that the above process can be aggregated into a function
	$u^{2,N}(\om, \theta_1, \theta_2)$ which is Borel measurable as a map from $\Om \x (\R_+)^2 \to \R$,
	\begin{equation} \label{eq:claim_aggregation}
		\mbox{uniformly continuous in}~\theta_1
		~\mbox{and}~ 
		u^{2,N}(\cdot, \tau_1, \tau_2) = Z^{2,N}_{\tau_1, \tau_2}, 
		~~\P_0-\mbox{a.s. for all}~\tau_1 \le \tau_2 \le N;
		~~~~
	\end{equation}
	and $u^{2,N}$ is increasing in $N$.

	Let
	\b*
		 v^{2,N}(\om, \theta_1, \theta_2)
		&:=&
		u^{2,N}(\om, \theta_1, \theta_2) 
		~+~
		\Phi_1(\om, \theta_1)
		~-~
		\lambda_1(\om_{\theta_1}),
	\e*
	and consider the optimal stopping problem
	\be \label{eq:stopping_pb_1}
		\sup_{\tau_1 \in \Tc^0_N} \E^{\P_0} \big[  v^{2,N}(\cdot, \tau_1, \tau_1) \big]
		&&\Big(=
			\sup_{\tau \in \Tc^a_N}
			\E^{\P_0} \big[ \Phi^{\lambda}(B, \tau) \big].
		\Big).
	\ee
	Denoted by $(Z^{1,N}_t)_{0 \le t \le N}$ the corresponding Snell envelop,
	which is $\F^a-$optional (or equivalently $\F^a-$predictable, since $\F^a$ is the augmented Brownian filtration),
	then $Z^{1,N}_t$ can be chosen to be $\F-$predictable (see e.g. Theorem IV.78 and Remark IV.74 of  Dellacherie \& Meyer \cite{DM}).
	Moreover, in view of \eqref{eq:Boundedness_Phi_lambda}, 
	by truncating it with $- C(1 + |\om_t| )$ from below  and with $C$ from above, we can assume that $Z^{1,N}_t$ is bounded between $- C(1 + |\om_t| )$ and $C$.
	Further, since $u^{2,N}$ is increasing in $N$, then for every $N_1 < N_2$, we know $Z^{1,N_2} \vee Z^{1,N_1}$ is still a Snell envelop of problem \eqref{eq:stopping_pb_1},
	then we can assume in addition and w.l.g. that $Z^{1,N}$ is increasing in $N$.
	Define $ v^{1,N}(\om, \theta_1) := Z^{1,N}(\om, \theta_1)$, 
	it follows that $v^{1,N}(\cdot), v^{2,N}(\cdot)$ are the required functionals.
	
	{\bf 2.} We now construct the measurable map $u^{2,N}$ satisfying \eqref{eq:claim_aggregation}.
	Let $\tau_1 \le \tau_2 \in \Tc^0_N$, define a random variable
	\be \label{eq:esssup_Z}
		\Zb^{2,N}_{\tau_1, \tau_2}
		&:=&
		 \esup_{\tau_3 \in \Tc^0_N, ~\tau_3 \ge \tau_2} 
		~\E^{\P_0} \Big[ \Phi_2 \big(B, \tau_1, \tau_3) - \lambda_2(B_{\tau_3}) ~ \Big| \Fc^a_{\tau_2} \Big].
	\ee
	Then, for every fixed $\tau_1$, $(\Zb^{2,N}_{\tau_1, \tau_2})_{\tau_2 \ge \tau_1}$ can be aggregated into a supermartingale, denoted by $Z^{2,N}_{\tau_1, t}$ (Lemma \ref{lemm:Snell}), such that $\Zb^{2,N}_{\tau_1, \tau_2} = Z^{2,N}_{\tau_1, \tau_2}$, $\P_0-$a.s. for each $\tau_2 \ge \tau_1$.
	Notice that $Z^{2,N}_{\tau_1, t}$ is $\F^a-$optional and equivalently $\F^a-$predictable, we can choose $Z^{2,N}_{\tau_1, t}$ to be $\F-$predictable (\cite[Theorem IV.78 and Remark IV.74]{DM}).
	Moreover, since $\Zb^{2,N}_{\tau_1, \tau_2}$ is increasing in $N$, $\P_0-$a.s.
	then for any $N_1 \le N_2$, $Z^{2,N_1} \vee Z^{2, N_2}$ is also an aggregated supermartingale for $(\Zb^{2,N_2}_{\tau_1, \tau_2})_{\tau_2 \ge \tau_1}$,
	and hence we can assume w.l.g. that $Z^{2,N}$ is increasing in $N$.
	Further, in view of \eqref{eq:Boundedness_Phi_lambda},
	by truncation, we can assume in addition and w.l.g. that
	$-C (1+ |\om_{t_1}| + |\om_{t_2}|) \le Z^{2,N}_{t_1, t_2} \le C$.
	
	Notice also that for two stopping times $\tau_1^1$ and $\tau_1^2$ smaller than $\tau_2$, we have
	\be \label{eq:Snell_envelop_coincid}
		\Zb^{2,N}_{\tau_1^1, \tau_2} 
		&=&
		\Zb^{2,N}_{\tau_1^2, \tau_2}, ~\P_0-\mbox{a.s.}
		~\mbox{on}~ A = \{\tau_1^1 = \tau_1^2 \}.
	\ee

	Further, since $\Phi_2(\om, \theta_1, \theta_2)$ is uniformly continuous in $\theta_1$, 
	denote by $\rho$ the continuity modulus.
	Then it follows by its definition in \eqref{eq:esssup_Z} that the family of random variables $\Zb^{2,N}_{\tau_1, \tau_2}$ is uniformly continuous w.r.t. 
	$\tau_1$, in sense that
	\b*
		\big| \Zb^{2,N}_{\tau^1_1, \tau_2} - \Zb^{2,N}_{\tau^2_1, \tau_2} \big|
		~\le~
		\rho(|\tau^1_1 - \tau^2_1 |),
		~\P_0 -
		\mbox{a.s. for stopping times}~ \tau_1^i \le \tau_2. 
	\e*
	We now define $u^{2,N}$ by
	\b*
		u^{2,N} (\om, \theta_1, \theta_2) &:=& Z^{2,N}_{\theta_1, \theta_2}(\om),
		~\mbox{for all}~ \theta_1 \in \Q,
	\e*
	and
	\b*
		u^{2,N}(\om, \theta_1, \theta_2) &:=& \limsup_{\Q \ni \theta_1' \to \theta_1} u^{2,N}(\om, \theta_1', \theta_2),
		~\mbox{for all}~ \theta_1 \notin \Q.
	\e*
	It is clear that $u^{2,N}$ is Borel measurable w.r.t. each variable since $Z^{2,N}_{\theta_1, \theta_2}(\om)$ is chosen to be $\F-$predictable.
	Furthermore, by \eqref{eq:Snell_envelop_coincid},
	{we have 
	$u^{2,N}(\om, \tau_1, \theta_2) = Z^{2,N}(\om, \tau_1, \theta_2)$ for all $\theta \ge \tau_1$, $\P_0-$a.s.,
	for every stopping times $\tau_1$ taking values in $\Q$.
	}
	Since we can approximate any stopping time by stopping times taking values in $\Q$, then by the uniform continuity of $\Zb^{2,N}_{\tau_1, \tau_2}$ w.r.t. $\tau_1$, we obtain that
	\b*
		\Zb^{2,N}_{\tau_1, \tau_2} ~=~ Z^{2,N}_{\tau_1, \tau_2} ~=~ u^{2,N}(\cdot, \tau_1, \tau_2)
		~\P_0-\mbox{a.s. for all stopping times}~ \tau_1 \le \tau_2 \in \Tc^1_N.
	\e*
	In particular, $u^{2,N}(\om, \theta_1, \theta_2)$ is uniformly continuous in $\theta_1$, $\P_0-$a.s.,
	which is the required functional in claim \eqref{eq:claim_aggregation}.
	\qed

	\begin{Remark}
		We notice that a general multiple optimal stopping problem has been studied in Kobylanski, Quenez \& Rouy-Mironescu \cite{Quenez},
		where the stopping times are not assumed to be ordered.
		In particular, they proved the existence of optimal multiple stopping times by a constructive method.
		Here we are in a specific context with Brownian motion and we are
		interested in finding a process $\Hb$ whose stochastic integral dominates the value process.
	\end{Remark}

\subsubsection{Proof of Theorem \ref{theo:main} \rmii  under Assumption \ref{Hyp:phi2} \rmii}
\label{subsubsec:duality2usc}

	Let $\Phi$ satisfy Assumption \ref{Hyp:phi} and Assumption \ref{Hyp:phi2} $\mathrm{(ii)}$, 
	i.e. $\omb \mapsto \Phi(\omb)$ is upper-semicontinuous and bounded from above.
	Define a metric $d$ of Polish space $\Omb$ by
	$$
		d(\omb, \omb') 
		~:=~ 
		\sum_{k=1}^n \big (|\theta_k - \theta_k'| 
		~+~
		\| \om_{\theta_k \wedge \cdot} - \om'_{\theta_k' \wedge \cdot} \| \big),
	$$
	and then define $\Phi_m:\Omb\to\R$ by
	\be\label{def:phim}
		\Phi_m(\omb) 
		:= \sup_{\omb'\in\Omb}\big\{\Phi(\omb')-m d(\omb,\omb')\big\}.~
	\ee
	Then $\Phi_m$ is a $d-$Lipschitz reward function, and satisfies in particular 
	Assumption \ref{Hyp:phi} and Assumption \ref{Hyp:phi2} $\mathrm{(i)}$. 
	Moreover, $\Phi_m(\omb)$ decreases to $\Phi(\omb)$ as $m$ goes to infinity for all $\omb\in\Omb$. 

	Denote by $P^m(\mu)$ and $D^m(\mu)$ the corresponding primal and dual values associated to the reward function $\Phi_m$.
	Since $\Phi_m$ satisfies Assumption \ref{Hyp:phi} and Assumption \ref{Hyp:phi2} $\mathrm{(i)}$, 
	we have proved in Section \ref{subsubsec:Duality2_1} the duality
	\b*
		P^m(\mu)
		&=&
		D^m(\mu).
	\e*
	Then by following the same line of argument as in Proposition \ref{prop:Phi_bounded}, we deduce that $P(\mu) = D(\mu)$.
	\qed

\appendix

\section{Appendix}

\subsection{On the canonical filtration on $\Omb$}

	We finally provide some properties of the canonical filtration $\Fbb = (\Fcb_t)_{t \ge 0}$ 
	of canonical space $\Omb$.
	Recall that the canonical element of $\Omb$ is denoted by $\big(B, T=(T_1, \cdots, T_n)\big)$,
	the $\sigma-$field $\Fcb_t$ is generated by the processes $B_{t\wedge \cdot}$ and $(T_k^t, k = 1, \cdots, n)$,
	where $T_k^t(\omb) := \theta_k \1_{\theta_k \le t} - \infty \1_{\theta_k > t}$ for all $\omb =\big(\om, \theta=(\theta_1, \cdots, \theta_n)\big) \in \Omb$.
	Equivalently, $\Fcb_t$ is generated by random variables $B_s$ and the sets $\{T_k \le s \}$ for all $k = 1, \cdots, n$ and $s \in [0,t]$.
	More importantly, $(T_k, k = 1, \cdots, n)$ are all $\Fbb-$stopping times.

	\begin{Lemma} \label{lemm:Fcb}
		The $\sigma-$field $\Fcb_{\infty}$ is the Borel $\sigma-$filed of $\Omb$.
		Moreover, the class of all bounded continuous, $\Fcb_t-$measurable functions on $\Omb$
		generates the $\sigma-$field $\Fcb_{t-} := \bigvee_{s < t} \Fcb_s$.
	\end{Lemma}
	\noindent{\bf Proof.} \rmi Since $T_k$ and $B$ are all $\Bc(\Omb)-$measurable, 
	one has $\Fcb_{\infty} \subseteq \Bc(\Omb)$.
	On the other hand, the process $(B_t, t \ge 0)$ generates the Borel $\sigma-$field $\Bc(\Om)$ and the collection of all sets $\{ T_k \le s\}$ generates the Borel $\sigma-$filed $\Bc(\Theta)$,
	it follows that $ \Bc(\Omb) = \Bc(\Om) \otimes \Bc(\Theta) \subseteq \Fcb_{\infty} $.
	
	\vspace{1mm}
	
	\noindent \rmii Let $t \ge 0$, denote $\Fc^B_t := \sigma( B_s, 0 \le s \le t)$, 
	$\Fc^{T_k}_t := \sigma \big( \{ T_k \le s\}, s \in [0,t] \big)$ and
	by $\Gc^{T_k}_t$ the $\sigma-$field generated by all bounded continuous and $\Fc^{T_k}_t-$measurable functions.
	First, for every $s <t$, it is clear that $\Fc^{T_k}_s \subset \Gc^{T_k}_t$, thus $\Fc^{T_k}_{t-} \subset \Gc^{T_k}_t$.
	Further, let $\phi: \R_+ \to \R$ be a bounded continuous function such that $\phi(T_k)$ is $\Fc^{T_k}_t-$measurable,
	then we have $\phi(t_1) = \phi(t_2)$ for every $t_1 \ge t_2 \ge t$.
	It follows that $\Phi(T_k)$ is $\Fc^{T_k}_{t-}-$measurable.
	Therefore, we have $\Fc^{T_k}_{t-} = \Gc^{T_k}_t$.
	Besides, it is well known that $\Fc^B_{t-} = \Fc^B_t$ is the $\sigma-$field generated by all bounded,  continuous and $\Fc^B_t-$measurable functions.
	It follows that $\Fcb_{t-} = \cup_{k=1}^n \Fc^{T_k}_{t-} \cup \Fc^B_{t-}$ is in fact the $\sigma-$field generated by all bounded, continuous and $\Fcb_t-$measurable functions.
	\qed

	\vspace{2mm}
	
	We now consider the filtration $\Fbb$.
	Let $t \ge 0$ and $\omb = (\om, \theta_1, \cdots, \theta_n) \in \Omb$, we introduce
	$[\omb]_t = (\om_{t \wedge \cdot}, [\theta_1]_t, \cdots, [\theta_n]_t)$ by
	$[\theta_k]_t := \theta_k \1_{\theta_k \le t} + \infty \1_{\theta_k > t}$.

	\begin{Lemma} \label{lemm:rcpd}
		\rmi $Y: \R_+ \x \Omb \to \R$ is $\Fbb-$optional if and only if it is $\Bc(\R_+ \x \Omb)-$measurable 
		and satisfies
		\be \label{eq:opt_rep}
			Y_s(\omb) &=& Y_s([\omb]_s) \mbox{ for all } s \ge 0 ~\mbox{and}~ \omb \in \Omb.
		\ee
		\rmii Consequently, $\Fcb_{T_k}$ is countably generated and every probability measure $\Pb$ 
		on $(\Omb, \Fcb_{\infty})$ admits a r.c.p.d. $(\Pb_{\omb})_{\omb \in \Omb}$ with respect to  $\Fcb_{T_k}$ which satisfies that
		\begin{itemize}
			\item[a)] $(\Pb_{\omb})_{\omb \in \Omb}$ is a family of conditional probabilities of $\Pb$ with respect to $\Fcb_{T_k}$,
			\item[b)] $\Pb_{\omb}( T_k= \theta_k, B_{T_k\wedge\cdot}= \om_{T_k\wedge\cdot}) = 1$ for all $\omb = (\om, \theta_1, \cdots, \theta_n) \in \Omb$.
		\end{itemize}
	\end{Lemma}
	\noindent{\bf Proof.}
	\rmi First, if $Y$ is $\Fbb-$optional, then $Y$ is measurable and $\Fbb-$adapted, i.e. $Y_s$ is $\Fcb_s-$measurable.
	Since $\Fcb_s$ is generated by $\omb \mapsto (\om_{s \wedge \cdot}, [\theta]_s)$, it follows that \eqref{eq:opt_rep} holds true.
	On the other hand, the process $(s, \omb) \mapsto (\om_{s \wedge \cdot}, [\theta]_s)$ is 
	adapted and c\`adl\`ag,
	and hence $\Fbb-$optional.
	Therefore, for every measurable process $\Yb$, the process $Y$ defined by \eqref{eq:opt_rep} is $\Fbb-$optional.
	
	\vspace{1mm}
	
	\noindent \rmii Notice that $\Bc(\Omb)$ is countably generated.
	And by the representation \eqref{eq:opt_rep}, the $\Fbb-$optional $\sigma-$field is generated by the map
	$(s, \omb) \in \R_+ \x \Omb \mapsto [\omb]_s \in \Omb$, and hence is also countably generated.
	Moreover, by Theorem IV-64 of Dellacherie \& Meyer \cite[p. 122]{DM}, we have
	\b*
		\Fcb_{T_k} &=& \sigma\{ B_{T_k \wedge \cdot}, ~ T_k \},
	\e*
	and hence $\Fcb_{T_k}$ is countably generated.
	Therefore, it follows by Theorem 1.1.6 in Stroock \& Varadhan \cite{SV} that every probability measure $\Pb$ on $(\Omb, \Fcb_{\infty})$ admits a r.c.p.d. with respect to  the $\sigma-$field $\Fcb_{T_k}$ satisfying the condition in item \rmii of the lemma.
	\qed

\subsection{Facts on the optimal stopping problem}

	 We next recall some useful results from the classical optimal stopping theory (see e.g. El Karoui \cite{EK},  Peskir \& Shiryaev \cite{PS}, Karatzas \& Shreve \cite{KS} etc.)
	 Let $(\Om^*, \Fc^*, \P^*)$ be an abstract complete probability space, which is equipped with a filtration $\F^* = (\Fc^*_t)_{t \ge 0}$ satisfy the usual conditions.
	Denote $\Fc^*_{\infty} := \vee_{t \ge 0} \Fc^*_t$ and by $\Tc^*$ the class of all $\F^*-$stopping times taking value in $[0, \infty)$.
	Let $Y$ be a $\F^*-$optional process defined on $\Om^*$ of class (D), i.e. the class $(Y_{\tau})_{\tau \in \Tc^*}$

	For every $\tau \in \Tc^*$, we denote by $\Tc^*_{\tau}$ the collection of all stopping times $\sigma$ in $\Tc^*$ such that $\sigma \ge \tau$.
	We then define a family of random variables
	$$
		Z^0_{\tau} ~:=~ \mbox{ess}\sup_{\sigma \in \Tc^*_{\tau}} \E\big[ Y_{\sigma} \big| \Fc_{\tau} \big],
		~~\mbox{for all}~\sigma \in \Tc^*_{\tau}.
	$$
	Then by the dynamic programming principle, 
	the family $(Z^0_{\tau})_{\tau \in \Tc^*}$ is a supermartingale system, i.e. 
	$Z^0_{\sigma} \ge  \E[Z^0_{\tau} | \Fc^*_{\sigma}]$ for all stopping times $\sigma \le \tau$ in $\Tc^*$.
	Using  Dellacherie \& Lenglart \cite[Thm. 6 and Rem. 7 c)]{DL}, it follows that one can find a l\`adl\`ag (left-limit and right-limit) optional process $Z = (Z_t)_{t \ge 0}$ which aggregates the family $(Z^0_{\tau})_{\tau \in \Tc^*}$, i.e.
	$$
		Z_{\tau} ~=~ Z^0_{\tau},~~\P^*-\mbox{a.s.}~~\mbox{for all}~ \tau \in \Tc^*.
	$$
	In particular, $Z = (Z_t)_{t \ge 0}$ is a strong supermartingale of class (D),
	and it is called the Snell envelope of process $Y$, or equivalently the minimum strong supermartingale dominating the optional process $Y$, 
	i.e. $Z_0 = \mbox{ess}\sup_{\tau \in \Tc^*} \E\big[ Y_{\tau} \big| \Fc_0 \big]$ and
	$Z_{\tau} \ge Y_{\tau}$ $\P^*$-a.s. for all $\tau \in \Tc^*$. 
	Using the optional cross-section theorem (see e.g. Theorem IV.86 in \cite{DM}), it follows that
	$$
		Z_t \ge Y_t,~~\mbox{for all}~t \ge 0, ~~\P^*-\mbox{a.s.}
	$$
	We summarize the above facts in the following lemma.
		
	\begin{Lemma} \label{lemm:Snell}
		Let $Y$ be an $\F^*-$optional process of class (D),
		then there is a  $\F^*-$optional l\`adl\`ag process $Z$, 
		which is the smallest strong supermartingale such that
		$Z_0 = \mbox{ess}\sup_{\tau \in \Tc^*} \E\big[ Y_{\tau} \big| \Fc_0 \big]$ and
		$Z_t \ge Y_t$ for all $t \ge 0$, $\P^*$-a.s. 
		In particular, one has $\E[Z_0] = \sup_{\tau \in \Tc_{\F^*}}  \E \big[Y_{\tau} \big]$.
	\end{Lemma}

	We next recall the Doob-Meyer decomposition for supermartingales without right continuity (see e.g. \cite[Theorem 20, Appendix I]{DellacherieMeyerB}
	or Mertens \cite[Theorem T3]{Mertens}).
	\begin{Lemma} \label{lemm:Mertens}
		Let $(\Om^*, \Fc^*, \P^*)$ be a probability space equipped with a filtration $\F^*= (\Fc^*_t)_{t \ge 0}$ satisfying the usual conditions,
		$X = (X_t)_{t \ge 0}$ be an $\F^*-$optional process class (DL) and an $\F^*-$supermartingale
		\footnote{ Here $X$ may not be of class (D), and $X$ is an $\F^*-$supermartingale if $\E[ X_{\tau} | \Fc_{\sigma}] \le X_{\sigma}$
		for all bounded $\F^*-$stopping times $\sigma \le \tau$.
		}.
		Then $X$ has a unique decomposition $X= X_0 + M - A$, where 
		$M_0 = A_0 = 0$, $M$ is a c\`adl\`ag $\F^*-$martingale, $A$ is an $\F^*-$predictable increasing process.
	\end{Lemma}

	The above decomposition allows one to define the quadratic co-variation of a (l\`adl\`ag) supermartingale with a continuous martingale in a pathwise way, as in Karandikar \cite{Karandikar}.
	Let us stay in the context of Lemma \ref{lemm:Mertens}, and assume that $W$ is a continuous martingale in the filtered probability space $(\Om^*, \Fc^*, \P^*, \F^*)$.
	We denote by $\F^{X,W} = (\Fc^{X,W}_t)_{t \ge 0}$ the raw filtration generated by $(X,W)$, i.e.
	$\Fc^{X,W}_t := \sigma( X_s, W_s~s \le t)$.
	Next, define 
	$$
		X^+_t ~:=~ \lim_{ \Q \ni s \searrow t} X_s
		= X_0 + M_t + A^+_t,
		~~~\mbox{with}~~A^+_t := \lim_{ \Q \ni s \searrow t} A_s.
	$$
	Then $X^+$ is clearly still an $\F^*-$supermartingale, and has c\`adl\`ag paths almost surely.
	Let $\tau^n_0 := 0$, $\tau^n_{i+1} := \inf\{ t \ge \tau^n_i~: |X^+_t - X^+_{\tau^n_i}| \ge 2^{-n} ~\mbox{or}~ |W_{\tau^n_i} -W_t| \ge 2^{-n} \}$,
	\be \label{eq:const_quad_cova}
		Q^n_t 
		~:=~
		\sum_{i=0}^{\infty} 
		\big(X^+_{\tau^n_{i+1} \wedge t} - X^+_{\tau^n_i \wedge t} \big)
		\big(W_{\tau^n_{i+1} \wedge t} - W_{\tau^n_i \wedge t} \big),
		~~~
		Q_t := \limsup_{n \to \infty} Q^n_t,
	\ee
	and finally $Q^-_0 := Q_0$ and $Q^-_t := \lim_{\Q \ni s \nearrow t} Q_s$ for $t > 0$.

	\begin{Lemma} \label{eq:pred_co_variat}
		In $(\Om^*, \Fc^*, \P^*, \F^*)$, the process $Q^-$ is indistinguishable from the quadratic co-variation $\langle X, W \rangle$ of $X$ and $W$ (or equivalent of $M$ and $W$).
		Moreover, $Q^-$ is an $\F^{X,W}-$predictable process.
	\end{Lemma}
	\proof Notice that the paths of $X^+$ is c\`adl\`ag, $\P^*-$a.s.,
	then following Theorem 3 of Karandikar \cite{Karandikar},
	the process $(Q_t)_{t \ge 0}$, taking values in $(-\infty, \infty]$, is indistinguishable from the quadratic co-variation between $X^+$ and $W$.
	Since $A^+$ has finite variation and $W$ is continuous, then $Q$ is also the quadratic co-variation between $X$ and $W$ (or between $M$ and $W$), and moreover $Q$ has continuous paths, $\P^*-$a.s.
	Then $Q^-$ and $Q$ are indistinguishable.
	
	Further, by its construction, it is clear that $Q^-$ is $\F^{X,W+}-$adapted, where $\F^{X,W+} =(\Fc^{X,W+}_t)_{t \ge 0}$ is the right-continuous filtration defined by $\Fc^{X,W+}_t := \lim_{s \searrow t} \Fc^{X,W}_s$.
	Since $Q^-$ is left-continuous, it follows that $Q^-$ is $\F^{X,W+}-$predictable,
	which is equivalent to be $\F^{X,W}-$predictable.
	\qed

	\vspace{2mm}

	We finally provide an equivalence result of the optimal stopping problems.
	Let $(\Om^*, \Fc^*, \P^*)$ be an abstract complete probability space, which is equipped with two filtrations $\F^* = (\Fc^*_t)_{t \ge 0}$ and $\G^* = (\Gc^*_t)_{t \ge 0}$, 
	where $\Fc^*_t \subseteq \Gc^*_t$ for every $t \ge 0$ and both filtrations satisfy the usual conditions.
	Denote $\Fc^*_{\infty} := \vee_{t \ge 0} \Fc^*_t$ and $\Gc^*_{\infty} := \vee_{t \ge 0} \Gc^*_t$.
	We denote further by $\Tc_{\F^*}$ the class of all $\F^*-$stopping times,
	and by $\Tc_{\G^*}$ the collection of all $\G^*-$stopping times.
	Let $Y$ be an $\F^*-$optional process defined on $\Om^*$ of class (D).

	\begin{Assumption} [K] \label{hyp:K}
		For every $t \ge 0$, every $\Gc_t-$measurable bounded random variable $X$ satisfies
		\b*
			\E \big[ X | \Fc^*_t \big]
			&=&
			\E \big[ X | \Fc^*_{\infty} \big],
			~\P^* -\mbox{a.s.}
		\e*
	\end{Assumption}

	\begin{Lemma} \label{lemm:OS_equiv}
		Under Assumption \ref{hyp:K} we have
		\b*
			\sup_{\tau \in \Tc_{\F^*}} \E[ Y_{\tau}] 
			&=&
			\sup_{\tau \in \Tc_{\G^*}} \E [ Y_{\tau}].
		\e*
	\end{Lemma}
	\noindent{\bf Proof.} The result follows by Theorem 5 of Szpirglas \& Mazziotto \cite{SzMa}.
	Notice that in \cite{SzMa}, $Y$ is assumed to be l\`adl\`ag,
	it can be easily generalized for any optional process by considering the Snell envelop of $Y$ w.r.t. the filtration $\F^*$, since its Snell envelop is l\`adl\`ag, $\P^*-$a.s.
	Further, $Y$ is also assumed to be positive in \cite{SzMa},
	which induces immediately the same result when $Y$ is of class (D) since 
	and process of class (D) can be dominated from below by a uniformly integrable martingale.
	\qed

\end{document}